\documentclass[11pt]{amsart}
\usepackage{amsmath,amssymb,txfonts}
\usepackage{amssymb}
\usepackage{amsmath}
\usepackage{mathrsfs}
\usepackage[shortlabels]{enumitem}
\usepackage{color}
\usepackage{extarrows}
\begin{document}
\newtheorem{theorem}{Theorem}    
\newtheorem{proposition}[theorem]{Proposition}
\newtheorem{conjecture}[theorem]{Conjecture}
\def\theconjecture{\unskip}
\newtheorem{corollary}[theorem]{Corollary}
\newtheorem{lemma}[theorem]{Lemma}
\newtheorem{sublemma}[theorem]{Sublemma}
\newtheorem{observation}[theorem]{Observation}
\theoremstyle{definition}
\newtheorem{definition}{Definition}
\newtheorem{notation}[definition]{Notation}
\newtheorem{remark}[definition]{Remark}
\newtheorem{question}[definition]{Question}
\newtheorem{questions}[definition]{Questions}
\newtheorem{example}[definition]{Example}
\newtheorem{problem}[definition]{Problem}
\newtheorem{exercise}[definition]{Exercise}

\numberwithin{theorem}{section} \numberwithin{definition}{section}
\numberwithin{equation}{section}

\def\earrow{{\mathbf e}}
\def\rarrow{{\mathbf r}}
\def\uarrow{{\mathbf u}}
\def\varrow{{\mathbf V}}
\def\tpar{T_{\rm par}}
\def\apar{A_{\rm par}}

\def\reals{{\mathbb R}}
\def\torus{{\mathbb T}}
\def\heis{{\mathbb H}}
\def\integers{{\mathbb Z}}
\def\naturals{{\mathbb N}}
\def\complex{{\mathbb C}\/}
\def\distance{\operatorname{distance}\,}
\def\support{\operatorname{support}\,}
\def\dist{\operatorname{dist}\,}
\def\Span{\operatorname{span}\,}
\def\degree{\operatorname{degree}\,}
\def\kernel{\operatorname{kernel}\,}
\def\dim{\operatorname{dim}\,}
\def\codim{\operatorname{codim}}
\def\trace{\operatorname{trace\,}}
\def\Span{\operatorname{span}\,}
\def\dimension{\operatorname{dimension}\,}
\def\codimension{\operatorname{codimension}\,}
\def\nullspace{\scriptk}
\def\kernel{\operatorname{Ker}}
\def\ZZ{ {\mathbb Z} }
\def\p{\partial}
\def\rp{{ ^{-1} }}
\def\Re{\operatorname{Re\,} }
\def\Im{\operatorname{Im\,} }
\def\ov{\overline}
\def\eps{\varepsilon}
\def\lt{L^2}
\def\diver{\operatorname{div}}
\def\curl{\operatorname{curl}}
\def\etta{\eta}
\newcommand{\norm}[1]{ \|  #1 \|}
\def\expect{\mathbb E}
\def\bull{$\bullet$\ }
\def\C{\mathbb{C}}
\def\R{\mathbb{R}}
\def\Rn{{\mathbb{R}^n}}
\def\Sn{{{S}^{n-1}}}
\def\M{\mathbb{M}}
\def\N{\mathbb{N}}
\def\Q{{\mathbb{Q}}}
\def\Z{\mathbb{Z}}
\def\F{\mathcal{F}}
\def\L{\mathcal{L}}
\def\S{\mathcal{S}}
\def\supp{\operatorname{supp}}
\def\pv{\operatorname{p.v.}}
\def\esslimsup{\operatorname{esslimsup}}
\def\essliminf{\operatorname{essliminf}}
\def\dist{\operatorname{dist}}
\def\essi{\operatornamewithlimits{ess\,inf}}
\def\esss{\operatornamewithlimits{ess\,sup}}
\def\xone{x_1}
\def\xtwo{x_2}
\def\xq{x_2+x_1^2}
\newcommand{\abr}[1]{ \langle  #1 \rangle}

\newcommand{\Norm}[1]{ \left\|  #1 \right\| }
\newcommand{\set}[1]{ \left\{ #1 \right\} }
\def\one{\mathbf 1}
\def\whole{\mathbf V}
\newcommand{\modulo}[2]{[#1]_{#2}}

\def\scriptf{{\mathcal F}}
\def\scriptg{{\mathcal G}}
\def\scriptm{{\mathcal M}}
\def\scriptb{{\mathcal B}}
\def\scriptc{{\mathcal C}}
\def\scriptt{{\mathcal T}}
\def\scripti{{\mathcal I}}
\def\scripte{{\mathcal E}}
\def\scriptv{{\mathcal V}}
\def\scriptw{{\mathcal W}}
\def\scriptu{{\mathcal U}}
\def\scriptS{{\mathcal S}}
\def\scripta{{\mathcal A}}
\def\scriptr{{\mathcal R}}
\def\scripto{{\mathcal O}}
\def\scripth{{\mathcal H}}
\def\scriptd{{\mathcal D}}
\def\scriptl{{\mathcal L}}
\def\scriptn{{\mathcal N}}
\def\scriptp{{\mathcal P}}
\def\scriptk{{\mathcal K}}
\def\frakv{{\mathfrak V}}

\def\scriptx{{\mathcal X}}
\def\scriptj{{\mathcal J}}
\def\scriptr{{\mathcal R}}
\def\scriptS{{\mathcal S}}
\def\scripta{{\mathcal A}}
\def\scriptk{{\mathcal K}}
\def\scriptp{{\mathcal P}}
\def\frakg{{\mathfrak g}}
\def\frakG{{\mathfrak G}}
\def\boldn{\mathbf N}

\title[Wiener Regularity]
{Wiener type regularity for  non-linear integro-differential equations}

\author{Shaoguang Shi}
\author{Guanglan Wang}
\author{ZhiChun Zhai}

\address{School of mathematics and statistics, Linyi University, Linyi 276005, China}
\email{shishaoguang@mail.bnu.edu.cn}
\address{School of mathematics and statistics, Linyi University, Linyi 276005, China}
\email{wangguanglan@lyu.edu.cn}
\address{Department of Mathematics and Statistics,  MacEwan University, Edmonton, Alberta T5J2P2, Canada}
\email{zhaiz2@macewan.ca}


\thanks{The authors were supported by Natural Science Foundation of China (\# 12271232, 12071197), Natural Science Foundation of Shandong Province (\# ZR2019YQ04, \#2020KJI002, \#ZR2021MA079). And they would like to thank Prof. Jie Xiao (Memorial University, Canada) for his suggestions and comments on this work.}

\date{\today}

\begin{abstract}
The primary purpose of this paper is to study the Wiener-type regularity criteria for non-linear equations driven by integro-differential operators, whose model is the fractional $p-$Laplace equation. In doing so, with the help of tools from potential analysis, such as fractional relative Sobolev capacities,  Wiener type integrals,   Wolff potentials,  $(\alpha,p)-$barriers, and  $(\alpha,p)-$balayages, we first prove the characterizations of the fractional thinness and the Perron boundary regularity. Then,  we establish a Wiener test and a generalized fractional Wiener criterion. Furthermore, we also prove the continuity of the fractional superharmonic function,  the fractional resolutivity, a connection between  $(\alpha,p)-$potentials and  $(\alpha,p)-$Perron solutions, and the existence of a capacitary function for an arbitrary condenser. \end{abstract}

\maketitle

\tableofcontents

\section{Introduction}\label{s1}
\subsection{Aim of the paper}
The theory of the boundary regularity of a domain,  initiated by Wiener \cite{Wi},  is fundamental and significant in potential analysis, which is connected closely with harmonic functions. Many researchers have extensively studied this topic. We mention a few of these works. In \cite{GZ}, Gariepy and Ziemer
established a regularity condition at the boundary of weak
solutions of the Dirichlet problem for quasilinear elliptic equations of second order in an open set $O\subset \mathbb{R}^n.$   In \cite{M2}, Maz\'{y}a established the sufficient part of the Wiener test for quasilinear elliptic equations. In \cite{KM}, Kilpel\"{a}inen and  Mal\'{y} proved the necessary part of the Wiener test for quasilinear elliptic equations.

In this paper, we focus on the
Wiener type regularity for the following non-linear integro-differential equation,
$$
\mathfrak{L}_{\alpha}u=0\quad\hbox{in}\quad \Omega\subset \mathbb R^n
$$
 and its related Dirichlet problems.
We assume throughout this paper that $\Omega\subset \mathbb{R}^{n}$ is an open domain unless otherwise specified. For  any $ \varphi\in C^\infty_c(\Omega)$ (the class of all infinitely differentiable functions with compact support in $\Omega$), the operator $\mathfrak{L}_{\alpha}$ can be defined as
$$
\langle\mathfrak{L}_{\alpha}, \varphi\rangle:=\int_{\mathbb{R}^{n}}\int_{\mathbb{R}^{n}}|u(x)-u(y)|^{p-2}(u(x)-u(y))(\varphi(x)-\varphi(y))K_{\alpha}(x,y)dxdy
$$
and the kernel function $K_{\alpha}: \mathbb R^n\times \mathbb R^n\longrightarrow \mathbb R$ is assumed to be measurable and satisfying
$$
\frac{1}{C|x-y|^{n+\alpha p}}\leq K_{\alpha}(x,y)\leq \frac{C}{|x-y|^{n+\alpha p}}\quad \hbox{for}\quad (\alpha,p)\in (0,1)\times (1,\infty)\quad \hbox{and}\quad C\geq 1.
$$
Much work has been carried out recently in studying $\mathfrak{L}_{\alpha}$ and the Dirichlet problems associated to $\mathfrak{L}_{\alpha}.$ Let us only review a few.

For the linear case, when $p = 2$ and when the kernel $K_\alpha$ equals to the Gagliardo kernel $K_\alpha(x, y) = |x - y|^{-(n+2\alpha)},$
 Caffarelli and his collaborators have studied the fractional obstacle problem in \cite{CF, CRS,CSS}.
In  \cite{KKL}, the authors studied the solutions defined via integration by parts with test functions, as viscosity
solutions or via comparison, and proved that the three notions coincide for bounded solutions.

 For the non-linear and possibly degenerate case,
 in \cite{CKP1, CKP2},
Castro et al. proved a general Harnack inequality and the regularity results for $\mathfrak{L}_{\alpha}u=0$ in $\Omega$ while $u=g$ in  $\Omega^c:=\mathbb{R}^{n}\setminus \Omega.$ In  \cite{KKP1,KKP2,KKP3,KKP4}, Korvenp\"{a}\"{a} et al.  systematically studied the corresponding obstacle problems. They  established  the existence and uniqueness of the solutions, the boundedness, continuity and H\"{o}lder continuity up to the boundary from the obstacle, the minimum of the corresponding
weak supersolutions becoming a weak supersolution, a comparison principle, a priori bounds, and the lower semicontinuity of supersolutions.
In \cite{KMS}, Kuusi et al. established the existence, regularity and potential theory for the non-linear integro-differential equations involving measure data. In \cite{P}, Palatucci studied  $(\alpha, p)$-superharmonic functions, and
the nonlocal counterpart of the Perron method in non-linear potential theory, and the connection among the fractional viscosity solutions, the weak solutions and the  $(\alpha, p)-$superharmonic functions.

Specially, when  $K_{\alpha}(x,y)=|x-y|^{-(n+\alpha p)}$, $\mathfrak{L}_{\alpha}$ becomes the fractional $p-$Laplace operator $(-\Delta_{p})^{\alpha}$, which can be understood as
$$
\langle(-\Delta_{p})^{\alpha}u, \varphi\rangle:=\int_{\mathbb{R}^{n}}\int_{\mathbb{R}^{n}}\frac{|u(x)-u(y)|^{p-2}(u(x)-u(y))(\varphi(x)-\varphi(y))}{|x-y|^{n+\alpha p}}dxdy,\quad \forall\varphi\in C_{c}^{\infty}(\mathbb R^n).
$$
Namely, $\mathfrak{L}_{\alpha}=(-\Delta_{p})^{\alpha}$ in the sense of distributions in this case. Alternatively, there exists a constant $C(n,\alpha,p)$ such that
$$
(-\Delta_{p})^{\alpha}u(x)=C(n,\alpha,p)\lim_{\epsilon\to 0}\int_{|y-x|\ge\epsilon}\frac{|u(x)-u(y)|^{p-2}(u(x)-u(y))}{|x-y|^{n+\alpha p}}\,dy.
$$
Then the classical $(\alpha, p)-$Laplace equation
$(-\Delta_{p})^{\alpha}u=0
$
is understood in the sense of
\begin{equation}\label{(1.1)}
\langle(-\Delta_{p})^{\alpha}u, \varphi\rangle =0\quad \forall\varphi\in C_{c}^{\infty}(\mathbb R^n)
\end{equation}
and the function $u$ is said to be a weak solution to $(\ref{(1.1)})$. Supersolution and subsolution of $(\ref{(1.1)})$ can be defined as
 $\langle(-\Delta_{p})^{\alpha}u, \varphi\rangle \geq0$ and $\langle(-\Delta_{p})^{\alpha}u, \varphi\rangle \leq0$ for nonnegative $\varphi\in C_{c}^{\infty}(\Omega),$
 respectively.
The related Dirichlet problem is often defined as
\begin{equation}\label{(1.2)}
\begin{cases}
(-\Delta_{p})^{\alpha}u=0\quad\hbox{in}\quad\Omega;\\
u=f\quad\hbox{in}\quad\Omega^{c}.
\end{cases}
\end{equation}
For simplicity, in this paper, we will study the fractional thinness and the fractional regularity of the integro-differential equations associated with $(-\Delta_{p})^{\alpha}.$ This means   $\mathfrak{L}_{\alpha}=(-\Delta_{p})^{\alpha}$ for the rest of this paper.
The operator $(-\Delta_{p})^{\alpha}$ can be viewed as a nonlocal version of the classical $p-$Laplace operator $-\Delta_{p}$ and many studies have been carried out on $(\ref{(1.1)})$ and $(\ref{(1.2)})$ for its wide applications in Physics, Biology and so on. To see this, we refer the interested readers to \cite{BL, CS1, CS2, CSS, CS, ILPS, IS, LL1, LL2, MRT, PXZ, S, TX, V, W, ZTZ} and the references therein.

Motivated by the above-mentioned excellent works on the regularity of non-linear integro-differential equations, this paper aims to establish the Wiener type regularity for non-linear integro-differential equations  using a newly developed theory of the fractional relative Sobolev capacity, which is  quite different from that of previously known results.  Before stating our main results, we review some basic definitions and preliminaries.

\subsection{Preliminaries}
The working space of this paper is the fractional Sobolev space
$W^{\alpha, p}(\Omega)$ on a domain $\Omega\subset\mathbb R^n$, which is defined as
$$
W^{\alpha, p}(\Omega)=L^{p}(\Omega)\cap \dot{W}^{\alpha, p}(\Omega)
$$
with $L^{p}(\Omega)$ the classical $p-$Lebesgue space on $\Omega$ and $\dot{W}^{\alpha, p}(\Omega)$ the space with the Gagliardo semi-norm
$$
[u]_{\dot{W}^{\alpha, p}(\Omega)}^{p}= \int_{\Omega}\int_{\Omega}\frac{|u(x)-u(y)|^{p}}{|x-y|^{n+\alpha p}}dxdy \quad\hbox{for all measurable function}\quad u\quad\hbox{on}\quad\Omega.
$$
$W_{0}^{\alpha,p}(\Omega)$ and $\dot{W}_{0}^{\alpha,p}(\Omega)$ are the
completion of $C^{\infty}_c(\Omega)$ under $\|\cdot\|_{W^{\alpha,p}(\Omega)}$ and $[\cdot]_{\dot{W}^{\alpha,p}(\Omega)}$, respectively.
\subsubsection{$(\alpha,p)$-harmonic function and $(\alpha,p)$-boundary regularity}\label{1.2}
 A function $u: \Omega\longrightarrow  \mathbb{R}$ is said to be $(\alpha, p)-$harmonic if $u\in C(\Omega)$ (all continuous functions on $\Omega$) is a solution to $(\ref{(1.1)})$ in the weak sense. As a concept closely related to supersolution, the $(\alpha, p)-$superharmonic function is given by a function $u$ satisfying the following three assumptions (see \cite{KKP4}):
$$\begin{cases}
{\rm(i)} \quad u \ \hbox{is lower semicontinuous (l.s.c)};\\
{\rm(ii)} \quad u\not\equiv \infty\ \hbox{in each component of}\  \Omega;\\
{\rm(iii)}\quad u\geq v \ \hbox{on} \ \partial O\longrightarrow  u\geq v \ \hbox{in}\ O\quad\hbox{for all open set}\ O\Subset \Omega\ \& \ \forall (\alpha, p)-\hbox{harmonic function}\ v\in C(\overline{O}).
\end{cases}
$$
Here, $O\Subset \Omega$ means $\overline{O}$ is a compact subset of $\Omega$.
More precisely, $(\alpha, p)-$superharmonic function is the viscosity supersolution for $(\ref{(1.1)})$(see \cite{KKP1} for more details).
A function $u$ is called $(\alpha, p)-$subharmonic if $-u$ is $(\alpha, p)-$superharmonic. For abbreviation, we use $H_{\alpha}(\Omega)$, $H_{\alpha}^{+}(\Omega)$ and $H_{\alpha}^{-}(\Omega)$ to represent the class of $(\alpha, p)-$harmonic functions, $(\alpha, p)-$superharmonic functions and $(\alpha, p)-$subharmonic functions in $\Omega$, respectively. It is immediately that $H_{\alpha}(\Omega)=H_{\alpha}^{+}(\Omega)\cap H_{\alpha}^{-}(\Omega)$.

In potential theory, a basic problem is the boundary regularity of $\Omega$. We recall that a point $x_{0}\in \partial\Omega$(the boundary of $\Omega$) is $(\alpha,p)-$regular if for every $u\in W^{\alpha,p}(\Omega)\cap C(\overline{\Omega})$,
there exists a function $f\in  H_{\alpha}(\Omega)$ with $$f-u\in W^{\alpha,p}_{0}(\Omega)\quad \hbox{and}
\lim_{x\longrightarrow  x_{0}}f(x)=u(x_{0}).$$
The existence and uniqueness of such  $f$ can be seen from \cite[Theorem 2.7]{Sh} and \cite[Theorem 1.1]{SX3}. $\Omega$ is called $(\alpha,p)-$regular if $x$ is $(\alpha,p)-$regular for each $x\in \partial\Omega$ and $(\alpha,p)-$irregular if $\Omega$ is not $(\alpha,p)-$regular.
A considerable amount of research has been performed  on the problem of fractional regularity during the last decade, such as Ros-Oton-Serra \cite{ROS} by developing a fractional analog of the Krylov boundary Harnack method, Lindgren-Lindqvist \cite{LL2} via Perron's method, Iannizzotto-Mosconib-Squassina \cite{IMS} utilizing barriers and Giacomoni-Kumar-Sreenadh \cite{GKS} using a suitable Caccioppoli inequality and the weak Harnack inequality. In this paper,
 the fractional regularity is studied by using a newly developed theory of the fractional relative Sobolev capacity.

Traditionally, regularity is defined in connection with Perron solutions. For the fractional case, the $(\alpha,p)-$Perron solution was first considered, for example, in
\cite{KKP4,LL2}. We  recall some related definitions as follows.
The upper $(\alpha,p)-$Perron solution $\overline{H}_{f}^{\alpha}$ and the lower $(\alpha,p)-$Perron solution $\underline{H}_{f}^{\alpha}$ of a function $f:\partial\Omega\rightarrow  [-\infty,+\infty]$ in $\Omega$ are given by
$$
\overline{H}_{f}^{\alpha}=\overline{H}_{f}^{\alpha}(\Omega)=\inf\left\{u: u\in \mathcal{U}_{f}^{\alpha}\right\},\quad\underline{H}_{f}^{\alpha}=\underline{H}_{f}^{\alpha}(\Omega)=\sup\left\{u: u\in \mathcal{L}_{f}^{\alpha}\right\},
$$
where the upper class $\mathcal{U}_{f}^{\alpha}$ and the lower class $\mathcal{L}_{f}^{\alpha}$ are defined as
$$
\mathcal{U}_{f}^{\alpha}=\left\{u: u\in H_{\alpha}^{+}(\Omega),\ u \,\,\hbox{is bounded below},\,\,\liminf_{x\longrightarrow  y}u(x)\geq f(y)\,\,\hbox{for all}\,\, y\in \partial \Omega\right\}
$$
and
$$
\mathcal{L}_{f}^{\alpha}=\left\{u: u\in H_{\alpha}^{-}(\Omega), \,\,u \,\,\hbox{is bounded above},\,\,\limsup_{x\longrightarrow  y}u(x)\leq f(y)\,\,\hbox{for all}\,\, y\in \partial \Omega\right\}.
$$
It follows that there hold the following three results:
\begin{itemize}	
\item[\rm (i)] $u\in \mathcal{U}_{f}^{\alpha}$ if and only if $-u\in \mathcal{L}_{f}^{\alpha}$,
\item[\rm (ii)] $\underline{H}_{f}^{\alpha}\leq \overline{H}_{f}^{\alpha}$,
\item[\rm (iii)]  $\overline{H}_{f}^{\alpha}\leq \overline{H}_{g}^{\alpha}$ if $f\leq g.$
\end{itemize}	

The $(\alpha,p)-$Perron's solution is an important tool to solve the Dirichlet problem (\ref{(1.2)}). It is a natural question to ask which one of the two $(\alpha,p)-$Perron solutions is the ``correct" solution to (\ref{(1.2)}).  It  follows from the comparison principle \cite[Theorem 16]{KKP4} that the $(\alpha,p)-$Perron solution coincides with the classical solution of (\ref{(1.2)}). That is, $\underline{H}_{f}^{\alpha}$ and $\overline{H}_{f}^{\alpha}$ are local solutions, see for example,  \cite[Theorem 22]{LL2} and also \cite[Theorem 2]{KKP4} which implies that $\underline{H}_{f}^{\alpha}$ and $\overline{H}_{f}^{\alpha}$ can be either identically $-\infty$ in $\Omega$, identically $+\infty$ in $\Omega$, or $(\alpha,p)-$harmonic in $\Omega$, respectively.

We say that a boundary point $x_{0}$ of $\Omega$ is Perron regular, if
$$
\lim_{x\longrightarrow  x_{0}}\overline{H}_{f}^{\alpha}(x)=f(x_{0})\quad \forall f\in C(\partial\Omega).
$$
$\Omega$ is called Perron regular if all points $x_{0}\in \partial\Omega$ are regular. Similarly, the same regularity is true if we replace $\overline{H}_{f}^{\alpha}$ with $\underline{H}_{f}^{\alpha}$ since $\overline{H}_{f}^{\alpha}=-\underline{H}_{-f}^{\alpha}$. We will show in Theorem \ref{Theorem 1.3} that the $(\alpha,p)-$regular boundary point agrees with the Perron regular boundary point if $\Omega$ is bounded.

\subsubsection{$(\alpha,p)$-Wiener type integral and $(\alpha,p)$-thinness}\label{1.3}
In the non-linear case, an important device in the central concepts of modern potential theory is the Wiener test (or Wiener criterion) introduced by Wiener \cite{Wi} to measure the boundary regularity in terms of capacity densities. We first adopt the following form of the Wiener type integral defined for an arbitrary set $E$ as
$$
\mathcal{W}_{p}^{\alpha}(E,x_{0})=\int_{0}^{1}\left(F_{\alpha, p}(x_{0},E,r)\right)^{\frac{1}{p-1}}\frac{dr}{r}:=\int_{0}^{1}\left(\frac{C_{\alpha, p}(E\cap B(x_{0},r), B(x_{0},2r))}{C_{\alpha, p}(B(x_{0},r), B(x_{0},2r))}\right)^{\frac{1}{p-1}}\frac{dr}{r},
$$
where $F_{\alpha, p}(x_{0},E,r)$ is the $(\alpha,p)$-capacity density function, $B(x_0,r)$ is the $x_0$-centered Euclidean ball with radius $r$ and $C_{\alpha, p}(E, \Omega)$ is the $(\alpha,p)$-capacity for any set $E\subset \Omega\subseteq \mathbb{R}^{n}$ which was defined as
$$
C_{\alpha, p}(E, \Omega)=\inf _{\text{open}\ O\supset E}C_{\alpha, p}(O, \Omega)=\inf _{\text{open}\ O\supset E}\sup_{\text{compact}\ K\subset O}C_{\alpha, p}(K, \Omega),
$$
where
$$ C_{\alpha, p}(K, \Omega):=\inf_{u\in \mathcal{X}^{\alpha, p}_{0}(K, \Omega)}[u]_{\dot{W}^{\alpha, p}(\Omega)}^{p}
\quad\hbox{
with}\quad
\mathcal{X}^{\alpha, p}_{0}(K, \Omega):= \left\{u: u\in \dot{W}^{\alpha, p}_{0}(\Omega)\quad\& \quad u\geq  \chi_{K}\right\}.
$$
Here $\chi_E$ stands for the characteristic function of a set $E$.
In practice, $\mathcal{X}^{\alpha, p}_{0}(K, \Omega)$ can be replaced by $$\mathcal{Y}^{\alpha, p}_{0}(K, \Omega):= \left\{u: u\in \dot{W}^{\alpha, p}_{0}(\Omega)\quad\& \quad 0\leq u\leq  1,\quad u=1\quad\hbox{on}\quad K\right\}$$ (see for example \cite{SX1}). Accordingly, a set $E$ is called $(\alpha,p)$-capacity zero, denoted by $C_{\alpha, p}(E)=0,$ if $C_{\alpha, p}(E\cap \Omega, \Omega)=0$ for all open sets $\Omega\subset \mathbb{R}^{n}$.   A property holds quasi everywhere (denoted by q.e. in the following) if it holds except for a set of  zero  $(\alpha,p)$-capacity. As we will also see from this paper (Theorem \ref{Theorem 1.2}-Theorem \ref{Theorem 1.6}) and the known results \cite{AH,AX,HKM2,LWXY,M,W1,X,X2,XY}, Sobolev type capacity is a powerful and useful tool in the study of potential theory, function spaces, harmonic analysis, partial differential equations and so on.

A set $E$ is
$$\begin{cases}
 (\alpha,p)-\hbox{thick at} \quad x_{0}\quad\hbox{if}\quad \mathcal{W}_{p}^{\alpha}(E,x_{0})=\infty,\quad\hbox{fractional thickness};\\
(\alpha,p)-\hbox{thin at} \quad x_{0}\quad\hbox{if}\quad\mathcal{W}_{p}^{\alpha}(E,x_{0})<+\infty,\quad\hbox{fractional thinness}.
\end{cases}
$$
The thinness was introduced by Adams and Meyer \cite{AM} to non-linear potential theory and then further studied in \cite{AH,HKM1,HW} to quasi-linear elliptic equations. We adopted in this paper to the fractional integro-differential equations, which partly inspired by some ideas from \cite{HKM2} for non-linear elliptic equations.

\subsubsection{$(\alpha,p)$-barrier and $(\alpha,p)$-balayage}\label{1.4}
A function $u$ is called fractional barrier (denoted by $(\alpha,p)$-barrier) \cite[Definition 25]{LL2}  relative to $\Omega$ at $x_{0}$ if
$$
u\in H_{\alpha}^{+}(\Omega);\quad
\liminf_{x\longrightarrow  y}u(x)>0\quad\hbox{for each}\quad y\in \partial\Omega\backslash \{x_{0}\};\quad
\lim_{x\longrightarrow  x_{0}}u(x)=0.
$$
By the minimum principle \cite[Lemma 2.5]{SX3}, an $(\alpha,p)$-barrier is always nonnegative and only $\mathbb{R}^{n}$ admits an $(\alpha,p)$-barrier that is not strictly positive. Furthermore, if the open set $O\subset \Omega$ and $u$ is a strictly positive $(\alpha,p)$-barrier relative to $\Omega$, then $u$ is an $(\alpha,p)$-barrier relative to $O$. One of our main results-Theorem \ref{Theorem 1.3} in this paper is to characterize the Perron regular boundary points in terms of the $(\alpha,p)$-barrier.

To introduce our main results, for a  function which is locally bounded from below, we  recall the $(\alpha,p)$-balayage from \cite{SX3} as
$$\widehat{\mathcal{B}}^{f}(x):=\widehat{\mathcal{B}}^{f}(\Omega)(x)=
\lim_{r\longrightarrow  0}\inf_{\Omega\cap B(x,r)}\mathcal{B}^{f}\big(\Omega\cap B(x,r)\big),
$$
where
$$\mathcal{B}^{f}:=\mathcal{B}^{f}(\Omega)=\inf \Psi_{f}:=\inf\left\{u: u\in H_{\alpha}^{+}(\Omega)\quad\hbox{and}\quad u\geq f\quad\hbox{in}\quad \Omega\right\}.$$
A relative version is defined for a nonnegative function $g$ on a  set $E\subset \Omega$ as $\widehat{\mathcal{B}}^{g}_{E}=\widehat{\mathcal{B}}^{f}$ (the $(\alpha,p)$-balayage relative to $E$) for
$$
f=\begin{cases}g \quad \hbox{on} \quad E;\\
0\quad \hbox{on}\quad\Omega\setminus E.
\end{cases}
$$
$\widehat{\mathcal{B}}^{1}_{E}$ is called the $(\alpha,p)$-potential of the set $E$ in $\Omega$. It follows from \cite[Theorem 1.2]{SX3} that $\widehat{\mathcal{B}}^{1}_{E}\in H_{\alpha}^{+}(\Omega)\cap H_{\alpha}(\Omega\backslash \overline{E}).$

Now, we are ready to state our main results.

\subsection{Statement of main results}
Our first result gives two equivalent characterizations of $(\alpha,p)-$thinness.
\begin{theorem}\label{Theorem 1.1}
Assume that $E\subset \mathbb{R}^{n}$ and $U(x_{0})$ denotes a neighborhood of $x_{0}$. Then the following three statements are equivalent
\begin{itemize}	
\item[\rm (i)] $E$ is $(\alpha,p)$-thin at $x_{0}$;

\item[\rm (ii)] There is a function $u\in H_{\alpha}^{+}(U(x_{0}))$ with $\liminf_{x\longrightarrow  x_{0}\ \&\ x\in E\backslash \{x_{0}\}} u(x)>u(x_{0});$

\item[\rm (iii)] There is a nonnegative function $u\in H_{\alpha}^{+}(U(x_{0}))$ such that $\widehat{\mathcal{B}}^{u}_{E\cap V(x_{0})}(x_{0})<u(x_{0})$  for $V(x_{0})\Subset U(x_{0})$.
\end{itemize}
\end{theorem}

The second result characterizes the fractional regularity when $\Omega$ is bounded.

\begin{theorem}\label{Theorem 1.3}
Assume that $C_{\alpha, p}(\{x_{0}\}, \Omega)=0$ for a finite point $x_{0}\in \partial\Omega$ and $\Omega$ is bounded. Then the following statements are equivalent
\begin{itemize}	
\item[\rm (i)] $x_{0}$ is Perron regular;

\item[\rm (ii)] There is an $(\alpha,p)-$barrier at $x_{0}$ relative to $\Omega$;

\item[\rm (iii)] $\widehat{\mathcal{B}}^{u}_{\overline{U}\backslash \Omega}(V)(x_{0})=u(x_{0})$ for nonnegative $u\in H_{\alpha}^{+}(V)$ and $U\Subset V$ are bounded open sets with $x_{0}\in U$;

\item[\rm (iv)] $\widehat{\mathcal{B}}^{1}_{\overline{B}\backslash \Omega}(2B)(x_{0})=1$ for all balls $B$ with $x_{0}\in B$;

\item[\rm (v)] $x_{0}$ is $(\alpha,p)-$regular.
\end{itemize}
\end{theorem}

If $\Omega$ is unbounded, then the point $\infty\in \partial\Omega$. Then all topological notions are therefore understood with respect to the space $\overline{\mathbb{R}^{n}}=\mathbb{R}^{n}\cup \{\infty\}$. It is a classical problem to  classify the Riemann surfaces or Riemannian manifolds,  which carry nonconstant bounded $(\alpha, p)-$superharmonic functions.
The following theorem shows that the existence of such functions is connected closely with the regularity of $\infty$ for the Dirichlet problem $(\ref{(1.2)})$.

\begin{theorem}\label{Theorem 1.4}
Assume that $\Omega$ is unbounded and the ball $B\subset \Omega$. Then the following statements are equivalent
\begin{itemize}	
\item[\rm (i)] $\infty$ is Perron regular for each $\Omega$;

\item[\rm (ii)] $\infty$ is Perron regular for ${\overline B}^{c}$;

\item[\rm (iii)] There is a nonconstant bounded function $u\in H_{\alpha}^{+}(\mathbb{R}^{n})$;

\item[\rm (iv)] $C_{\alpha, p}(B, \mathbb{R}^{n})>0$ for each ball $B\Subset \Omega$;

\item[\rm (v)] $C_{\alpha, p}(B, \mathbb{R}^{n})>0$ for some ball $B\Subset \Omega$.
\end{itemize}
\end{theorem}

 With the help of the local nature of $(\alpha,p)-$superhamonic functions, the continuity of $(\alpha,p)-$balayage and the existence of the solution to (\ref{(1.2)}) with Sobolev boundary values that will be established in section \ref{s5},   we derive  the following fractional Wiener test  from Theorem \ref{Theorem 1.1}.

\begin{theorem}\label{Theorem 1.2}
Let $x_{0}\in \partial\Omega$ be a finite boundary point. Then
$$
x_{0} \, \hbox{is}\,(\alpha,p)-\hbox{regular} \Longleftrightarrow \mathcal{W}_{p}^{\alpha}(\Omega^{c}, x_{0})=+\infty \Longleftrightarrow \Omega^{c}\, \hbox{is}\, (\alpha,p)-\hbox{thick at}\, x_{0}.
$$
\end{theorem}

The last main results is a generalization of the classical Wiener Criterion, which need two basic notions.
Write
 $$(\alpha,p)-\esss f_{B(x_{0},r)}=\inf \left\{t:f\leq t\quad  q.e.\, \hbox{in}\quad B(x_{0},r)\right\},\,\,\,\,\overline{f(x_{0})}=\inf_{r>0}\left((\alpha,p)-\esss f_{B(x_{0},r)}\right)$$
and
$$
F_{\varepsilon}=\left\{x: f(x)\geq \overline{f(x_{0})}-\varepsilon\quad\hbox{for}\quad\varepsilon>0\right\}.
$$
Then,  the point $x_{0}$ is called an $(\alpha,p)-$Wiener point of $f$ if $F_{\varepsilon}$ is not $(\alpha,p)-$thin at $x_{0}$.

Denote by
$$
\Psi_{f,u_{0}}(\Omega)=\left\{u\in W^{\alpha,p}(\Omega):u\geq f\quad \hbox{in}\quad \Omega,\quad  u_{0}\in W^{\alpha,p}(\Omega)\quad \hbox{with}\quad u-u_{0}\in W_{0}^{\alpha,p}(\Omega)\right\}.
$$
The obstacle problem with obstacle $f$ and boundary value $u_{0}$ for $(\ref{(1.1)})$(denoted by $\Phi_{f,u_{0}}(\Omega)$) is to find a function $u\in \Psi_{f,u_{0}}(\Omega)$ such that
$$
\langle(-\Delta_{p})^{\alpha}u, \varphi\rangle\geq 0\quad\hbox{for}\quad\varphi\in \Psi_{f,u_{0}}(\Omega).
$$
Such a function $u$ is called the solution to $\Phi_{f,u_{0}}(\Omega).$
Let $f$ be bounded. A function $u$ is said to be a local solution to the obstacle problem at the point $x_{0},$ denoted by $\Phi^{x_0}_{f,u_{0}}(\Omega),$  if there is an open neighborhood $\Omega$ of $x_{0}$ such that
$$u\in W^{\alpha,p}(\Omega),\quad u\geq f\quad  q.e.$$
and
$$
\langle(-\Delta_{p})^{\alpha}u, \varphi\rangle\geq 0\quad \hbox{for}\quad \varphi\in W_{0}^{\alpha,p}(\Omega) \quad \hbox{with}\quad  u+\varphi\geq f\quad  q.e..
$$
 For more information about the obstacle problem, see, for example,  \cite{KKP2} and \cite{Sh}.

\begin{theorem}\label{Theorem 5.1} The following two statements hold.

\begin{itemize}	
\item[\rm (i)] If $x_{0}$ is an $(\alpha,p)-$Wiener point of $f,$ then each local solution to the obstacle problem $\Phi^{x_0}_{f,u_{0}}(\Omega)$ is continuous at $x_{0}.$

\item[\rm (ii)] If $x_{0}$  is not an $(\alpha,p)-$Wiener point of $f$, then there exits a local solution of $\Phi^{x_0}_{f,u_{0}}(\Omega)$ which can not be  continuous at $x_{0}$.
\end{itemize}
\end{theorem}

\subsection{Plan of the paper}
The rest of this paper is organized as follows. In Section \ref{S2}, we  provide the proof of  our main results.
More specifically, in section \ref{s2},  based on one lemma concerning the $(\alpha,p)-$thinnes/ thickness and another one discussing when $\widehat{\mathcal{R}}^{f}=\widehat{\mathcal{B}}^{f},$ we will prove Theorem \ref{Theorem 1.1} which provides  some characterizations for $(\alpha,p)$-thinness and $(\alpha,p)$-boundary regularity.
In Section \ref{s3}, we establish Theorem \ref{Theorem 1.3} which provides regularity characterizations when $\Omega$ is bounded. Moreover, as a byproduct, we  show that the regularity is a local property. When $\Omega$ is unbounded,  in Section \ref{s4}, after establishing the weak compactness in $\dot{W}^{\alpha,p}_0(\Omega),$ we prove Theorem \ref{Theorem 1.4} and its corollary. In Section
\ref{s5}, based on four technical lemmas which are significant for the understanding of the $(\alpha,p)-$supperhamonic functions, the Perrson regularity, the $(\alpha,p)-$balayage, and the solution of (\ref{(1.2)}) with Sobolev boundary values,  we prove Theorem \ref{Theorem 1.2} which provides the fractional Wiener test. In Section \ref{sec2.5}, we derive Theorem \ref{Theorem 5.1} from Theorem \ref{Theorem 1.1}.
 To see the significance of our main results and technical lemmas,
in Section \ref{s8} we derive more regularity conditions.  After proving a generalized comparison lemma, in Section \ref{s6}, we prove Theorem \ref{Theorem 1.5} concerning the continuity of $(\alpha,p)-$supperharmonic functions. Based on  a technical   lemma about the  uniform convergence for upper $(\alpha,p)-$Perron solutions, we prove the resolutivity( Theorem \ref{Theorem 1.6}) in Section \ref{s7}.  In the rest of Section \ref{s8}, we prove
a connection between  $(\alpha,p)-$potentials and  $(\alpha,p)-$Perron solutions,
and the existence of a capacitary function for an arbitrary condenser
using Theorem \ref{Theorem 1.1},  Theorem \ref{Theorem 1.5},
Lemma \ref{lemma 3.4}, Lemma \ref{lemma 3.6}
 and Lemma \ref{lemma 5.1}.

In the forthcoming discussions, $A\lesssim B$ ($A\gtrsim B$) means $A\leqslant CB$ ($A\geqslant CB$) for a positive constant $C$ which may change from line to line and $A\thickapprox B$ amounts to $A\lesssim B\lesssim A$. $u^{+}=\max\{u, 0\}$, $u^{-}=\min\{u, 0\}.$

\section{Proof of main results}\label{S2}
This section is devoted to the proof of our main results.   We begin with the proof of Theorem  \ref{Theorem 1.1}.
\subsection{Proof of Theorem  \ref{Theorem 1.1}:  characterizations of fractional thinness }\label{s2}

 We only need to show the cycle $(i)\Longrightarrow  (ii)\Longrightarrow  (iii)\Longrightarrow  (i).$  To do so, some technical lemmas are needed.
\subsubsection{Important Lemmas}
\begin{lemma}\label{lemma 2.1}
Suppose that $E$ is a Borel set, then
$$
\begin{cases}
E\,\hbox{is}\, (\alpha,p)-\hbox{thin at}\,\, x_{0}\Longrightarrow \hbox{there exists an open set}\, O\supset E\backslash \{x_{0}\}\,\hbox{and}\, O\, \hbox{is}\,\,(\alpha,p)-\hbox{thin at}\, x_{0};\\
E\, \hbox{is}\, (\alpha,p)-\hbox{thick at}\, x_{0}\Longrightarrow \hbox{there exists a compact set}\, K\subset E\backslash \{x_{0}\}\, \hbox{and}\, K\, \hbox{is}\, (\alpha,p)-\hbox{thick at}\, x_{0}.
\end{cases}
$$
\end{lemma}

\begin{proof}
Let $B_{i}=B(x_{0},2^{-i})$ for $i=1,2,\cdots, $ and pick an open set $O_{i}\subset B_{i}$ with $E\cap B_{i}\subset O_{i}$ and
$$
\left(F_{\alpha,p}(x_{0},E,2^{-i})\right)^{\frac{1}{p-1}}\geq \left(\frac{C_{\alpha,p}(O_{i}, B_{i-1})}{C_{\alpha,p}(B_{i}, B_{i-1})}\right)^{\frac{1}{p-1}}-2^{-i},
$$
where  $F_{\alpha, p}(x_{0},E,r)$ is the $(\alpha,p)$-capacity density function.
Without loss of generality, we   assume that the sequence $\{O_{i}\}$ is decreasing and that $E\subset B_{1}$. Therefore, $E\backslash \{x_{0}\}\subset O=\cup_{i}(O_{i}\backslash \overline{B}_{i+2})$. Moreover, $C_{\alpha,p}(O\cap B_{1}, B_{0})\leq C_{\alpha,p}(B_{1}, B_{0})$, and hence
$$
C_{\alpha,p}(O\cap B_{i}, B_{i-1})\leq C_{\alpha,p}(O_{i-1}\cap B_{i}, B_{i-1})\leq CC_{\alpha,p}(O_{i-1}, B_{i-2}).
$$
Accordingly,
$$
\sum_{i=1}^{\infty}\left(F_{\alpha,p}(x_{0},O,2^{-i})\right)^{\frac{1}{p-1}}\leq C+C\sum_{i=2}^{\infty}\left(\left(F_{\alpha,p}(x_{0},E,2^{-i})\right)^{\frac{1}{p-1}}+2^{-i}\right)<+\infty,
$$
which is the desired result by \cite[Lemma 2.10]{SX3}.

Similarly, we can prove  the second assertion by recalling $C_{\alpha,p}(E, \Omega)=\sup_{K\subset E}C_{\alpha,p}(K, \Omega)$ with $K$ being compact.
\end{proof}

In order to state the next lemma, we need  another version of $(\alpha,p)$-balayage. For a function $f$, which is locally bounded from below, let
$$
\mathcal{R}^{f}=\mathcal{R}^{f}_{\Omega}=\inf \left\{u: u\in H_{\alpha}^{+}(\Omega)\quad \hbox{and}\quad u\geq f\quad  q.e.\quad in \quad \Omega\right\}.
$$
Then the quasi $(\alpha,p)$-balayage can be defined as
$$
\widehat{\mathcal{R}}^{f}(x)=\widehat{\mathcal{R}}^{f}_{\Omega}(x)=\liminf_{y\longrightarrow  x}\mathcal{R}^{f}_{\Omega}(y).
$$
If $f$ is bounded, then $\widehat{\mathcal{R}}^{f}\in H_{\alpha}^{+}(\Omega)$ by a slight modification of \cite[Lemma 2.7]{SX3} and $\widehat{\mathcal{R}}^{f}_{\Omega}\geq f$ q.e. in $\Omega$ by \cite[Theroem 1.3]{SX3}.
$\widehat{\mathcal{R}}^{f}_{\Omega}$ is the smallest function in $H_{\alpha}^{+}(\Omega)$ above the obstacle $f$ in $\Omega$ q.e.. If $\Omega$ is bounded, then $f$ is bounded, and hence $\widehat{\mathcal{R}}^{f}_{\Omega}\in W_{\alpha,p}(\Omega)$ by \cite[Proposition 2.11]{Sh} and \cite[Theorems 11 \& 13]{KKP4}.

By recalling the definition of $(\alpha,p)-$balayage $\widehat{\mathcal{B}}^{f}$ in Section \ref{s1}, it is interesting to know whether sets of $(\alpha,p)-$capacity zero can be neglected. Namely, whether
$$
\widehat{\mathcal{R}}^{f}=\widehat{\mathcal{B}}^{f}?
$$
Below is a desired solution to this question.
\begin{lemma}\label{lemma 2.2}
Let $\Omega$ be a bounded open set. Then the following statements hold.
\begin{itemize}
\item[\rm (i)] If $\widehat{\mathcal{R}}^{f}\in \dot{W}^{\alpha,p}(\Omega),$ then $\widehat{\mathcal{R}}^{f}=\widehat{\mathcal{B}}^{f}$;

\item[\rm (ii)]If  $f$ is nonnegative with $\supp f\Subset \Omega$, then $ \widehat{\mathcal{R}}^{f}=\widehat{\mathcal{B}}^{f}$;

\item[\rm (iii)] $\widehat{\mathcal{R}}^{f}$ is a local solution of $\Phi^{x_0}_{f,u_{0}}(\Omega)$.
\end{itemize}
\end{lemma}
\begin{proof}
For abbreviation, we write $\widehat{\mathcal{R}}, \widehat{\mathcal{B}}$ instead of $\widehat{\mathcal{R}}^{f}$ and $\widehat{\mathcal{B}}^{f}$, respectively. It remains to prove that $\widehat{\mathcal{B}}\leq \widehat{\mathcal{R}}$ since $\widehat{\mathcal{B}}\geq \widehat{\mathcal{R}}$ q.e. is immediately.

By letting $S=\left\{x\in \Omega: \widehat{\mathcal{R}}(x)<f(x)\right\}$, we get from the fundamental convergence theorem \cite[Theorem 1.3]{SX3} that $C_{\alpha,p}(S, \Omega)=0.$ Hence,  there exists a nonnegative l.s.c function $u\in W_{\alpha,p}(\mathbb{R}^{n})$ with $u(x)=\infty$ for $x\in S$.

Let $f_{i}=\widehat{\mathcal{R}}+\frac{1}{i}u$ and  $v_{i}$ be the solution to $\Phi_{f_{i}, f_{i}}(\Omega)$. Then \cite[Theorem 2.5\,,\,Theorem 3.1 \& Corollary 3.7]{Sh} allow us to assume that $v_{i}\in H_{\alpha}^{+}(\Omega).$ Hence,  since $f_{i}$ is l.s.c., we get
$$
v_{i}(x)=\essliminf_{y\longrightarrow  x}v_{i}(y)\geq \essliminf_{y\longrightarrow  x}f_{i}(y)\geq f_{i}(x)\geq f(x)\quad\forall x\in \Omega
$$
according to  \cite[Proposition 2.11]{Sh}.  Accordingly, $\widehat{\mathcal{B}}\leq \lim v_{i}=v.$
Moreover, the fact $f_{i}\longrightarrow  \widehat{R}$ in $\dot{W}_{\alpha,p}(\Omega)$, the uniqueness of solution to $\Phi_{f_{i},f_{i}}$ and \cite[Proposition 2.15]{Sh} imply that $v=\widehat{\mathcal{R}}$ q.e.. Finally, according to
\cite[Corollary 3.9]{Sh},  we get $$\widehat{\mathcal{R}}(x)=\essliminf_{y\longrightarrow  x}v(y)\geq \essliminf_{y\longrightarrow  x}\widehat{\mathcal{B}}(y)=\widehat{\mathcal{B}}(x)\quad  \forall x\in \Omega,$$  which is (i).

To prove (ii), we only need to show $\widehat{\mathcal{R}}\in \dot{W}_{\alpha,p}(\Omega)$ according to  (i). By choosing a neighborhood $U$ of $\partial \Omega$ such that
$$
\begin{cases}
\overline{U}\cap \supp f=\phi;\\
\mathbb{R}^{n}\backslash U \, \hbox{is not}\, (\alpha,p)-\hbox{thin at}\, \hbox{every}\quad  x\in \partial U,
\end{cases}
$$
it follows from \cite[Lemma 2.7]{SX3} that $\widehat{\mathcal{R}}\in H_{\alpha}(\Omega\backslash \supp f)$ and hence $\widehat{\mathcal{R}}\in H_{\alpha}(U\cap \Omega)$.

Let
$$\left\{
\begin{array}{ll}
\partial U\cap\Omega\subset \Omega_{i}\Subset \Omega_{i+1}\Subset \Omega\quad \hbox{with}\quad \cup_{i}\Omega_{i}=\Omega;\\
\varphi= \widehat{\mathcal{R}}\,\,\,\, \hbox{in}\,\,\,\, \partial U\cap\Omega\quad \hbox{with}\quad \supp \varphi\Subset \Omega\quad  \hbox{for}\quad \varphi\in C(\Omega)\cap \dot{W}^{\alpha,p}(\Omega);\\
 h_{i}\in H_{\alpha}(U\cap \Omega_{i})\quad \hbox{with}\quad  h_{i}-\varphi\in  \dot{W}_{0}^{\alpha,p}(U\cap \Omega_{i}).
\end{array}\right.
$$
Then the comparison principle \cite[Lemma 6]{KKP4}(see also \cite[Theorem 15]{P}) gives
$$
\widehat{\mathcal{R}}(x)\geq h_{i+1}(x)\geq h_{i}(x)\quad\hbox{for}\quad x\in U\cap \Omega_{i}
$$
since $\widehat{\mathcal{R}}\in W^{\alpha,p}(U\cap \Omega_{i})$\cite[Theorem 1, Lemma 5 \& Theorem 11]{KKP4} with  $$\widehat{\mathcal{R}}\geq 0 \quad \hbox{and}\quad \widehat{\mathcal{R}}=\varphi=h_{i} \quad \hbox{in}\quad  \partial U\cap \Omega_{i}=\partial U\cap \Omega.$$ Consequently,  \cite[Theorem 15]{KKP4} implies $\lim_{i}h_{i}=h\in H_{\alpha}(U\cap \Omega).$   Furthermore, we get $$[h_{i}]_{\dot{W}_{\alpha,p}(U\cap \Omega_{i})}^{p}\leq C[\varphi]_{\dot{W}_{\alpha,p}(\Omega)}^{p},$$ which shows $h\in \dot{W}^{\alpha,p}(U\cap \Omega)$.

Next, we claim that $h=\widehat{\mathcal{R}}$ in $U\cap \Omega$. In fact, $h\leq \widehat{\mathcal{R}}$ is immediately. To prove the inverse, let
$$
\widetilde{h}=\begin{cases}
\min\{h,\widehat{\mathcal{R}}\}\,\, \hbox{in}\,\,  U\cap\Omega;\\
\widehat{\mathcal{R}}\,\,  \hbox{in}\,\, \Omega\backslash U.
\end{cases}
$$
Then $\widetilde{h}$ is l.s.c in $\Omega$ by the fact $\lim_{y\longrightarrow  x}h(y)=\widehat{\mathcal{R}}(x)$ for all  $ x\in \partial U\cap \Omega.$ It follows from  \cite[Proposition 2.21]{Sh}  that $\widetilde{h}\in H_{\alpha}^{+}(\Omega).$ Whence we get $\widetilde{h}\geq \widehat{\mathcal{R}}$ in $\Omega$, and finally $h= \widehat{\mathcal{R}}$ in $U\cap \Omega$, which complete the proof of (ii).

(iii). Let $x_{0}\in \Omega$ and $B=B(x_{0},r)\Subset \Omega$. Then $\widehat{\mathcal{R}}^{f}\in W^{\alpha,p}(B)$ by \cite{KKP4}. Suppose that $u$ is the solution to $\Phi_{f, \widehat{\mathcal{R}}}(B)$. Then  it follows from  \cite[Theorem 3.1 \& Corollary 3.7]{Sh} that $u\in H_{\alpha}^{+}(B).$ Hence we obtain $\widehat{\mathcal{R}}\geq u$ in $B$  according to \cite[Lemma 2.6]{Sh} since $\min\left\{u, \widehat{\mathcal{R}}\right\}-\widehat{\mathcal{R}}\in W_{0}^{\alpha,p}(B)$.

Next, we need to show $u\geq \widehat{\mathcal{R}}$ in $B$ in order to  complete our proof. Denote by $\{\varphi_{i}\}$ an  increasing  sequence with
$\{\varphi_{i}\}\subset C_{c}^{\infty}(\mathbb{R}^{n})$ and $\varphi_{i}\longrightarrow  \widehat{\mathcal{R}}$ in $\partial B.$ Let
$
h_{i}\in C(\overline{B})\cap H_{\alpha}(B)$ be the unique function such that $ h_{i} =\varphi_{i}$ in $\partial B.
$
It follows from  \cite[Lemma 2.8]{Sh}  that  $\min\{u-h_{i},0\}\in W_{0}^{\alpha,p}(B)$ and hence $u\geq h_{i}$ in $B.$  This  implies that $\liminf_{y\longrightarrow  x}u(y)\geq \widehat{\mathcal{R}}(x)$ for all $x\in \partial B.$ Accordingly, the function
$$
v=\begin{cases}
\min\{\widehat{\mathcal{R}},u\}\quad\hbox{in}\quad B;\\
\widehat{\mathcal{R}}\quad\hbox{in}\quad\Omega\backslash B
\end{cases}
$$
is l.s.c and $v\in H_{\alpha}^{+}(\Omega)$ due to \cite[Proposition 2.21]{Sh}. Consequently, $u\geq v\geq \widehat{\mathcal{R}}$ in $B$ as desired.
\end{proof}

Having proved the previous two critical lemmas, we are ready to prove Theorem \ref{Theorem 1.1}.
\subsubsection{Proof of Theorem \ref{Theorem 1.1}}
(i) $\Longrightarrow$ (ii). Suppose that $E$ is $(\alpha,p)-$thin at $x_{0}\notin E.$ Then, we can assume that $E$ is open according to Lemma \ref{lemma 2.1}. Taking
$$
\begin{cases}
B_{i}=B(x_{0},r_{i})\quad\hbox{with}\quad r_{i}=2^{-i};\\
E_{i}=E\cap B_{i};\\
u=\widehat{\mathcal{B}}^{1}_{E_{j}}(B_{j-2})\quad\hbox{with}\quad(-\Delta_{p})^{\alpha}u=\mu\quad\hbox{and}\quad j\geq 2\quad\hbox{be an integer},
\end{cases}
$$
we conclude that $u\geq 1$ on $E_{j}$. In order to obtain (ii), it suffices to show that $u(x_{0})<1$. It follows from  \cite[Lemma 2.8]{SX3} and \cite[Theorem 2.2]{SX1}  that
$$
\left(\inf_{B_{j}}u\right)^{p-1}r_{j}^{n-\alpha p}\leq C\left(\inf_{B_{j}}u\right)^{p-1}C_{\alpha,p}\left(\left\{u>\inf_{B_{j}}u\right\}, B_{j-2}\right)\leq C\mu(B_{j-2})\leq C\mu(B_{j-1})
$$
and consequently,
$$
\inf_{B_{j}}u\leq C\left(\frac{\mu(B_{j-1})}{r_{j-1}^{n-\alpha p}}\right)^{\frac{1}{p-1}}.
$$
On the other hand, \cite[Lemma 5.1]{SX2} shows that,  for $i>j-2,$ $$\mu(B_{i})\leq C C_{\alpha,p}(E_{i}, B_{j-2})\leq C C_{\alpha,p}(E_{i}, B_{i-1}).$$ The previous inequality  together with
\cite[Theorem 1.2]{KMS} and the fact $u\geq 1$ implies that there exists $C:=C(n,\alpha,p)>0$ such that
$$
u(x_{0})\leq C\inf_{B_{j}}u+C\textbf{W}_{\alpha,p}^{\mu}(x_{0}, r_{j-1})\leq C\sum_{i=j-1}^{\infty}\left(\frac{C_{\alpha,p}(E_{i}, B_{i-1})}{r_{i}^{n-\alpha p}}\right)^{\frac{1}{p-1}}\leq\frac{1}{2}
$$
as desired by taking $j$ large enough.

(ii) $\Longrightarrow$ (iii). We give the proof only for the case $x\in \overline{E\backslash \{x_{0}\}}$ since the case $x\notin \overline{E\backslash \{x_{0}\}}$
is straightforward. We can assume that $\liminf_{x\rightarrow  x_{0}, x\in E\backslash \{x_{0}\}}u(x)>1>u(x_{0})$ for the $u$ given in (ii). By choosing an open set $V(x_{0})\Subset U(x_{0})$ with $u|_{(E\cap V(x_{0}))\backslash \{x_{0}\}}>1$, Lemma \ref{lemma 2.2} implies
$$
\widehat{\mathcal{B}}^{1}_{E\cap V(x_{0})}(x_{0})=\widehat{\mathcal{B}}^{1}_{(E\cap V(x_{0}))\backslash \{x_{0}\}}(x_{0})\leq u(x_{0})<1,
$$
which is (iii).

(iii) $\Longrightarrow$ (i). Without loss of generality, we  assume that $U(x_{0})$ is bounded and $(\alpha, p)-$regular. The choquet topological lemma \cite[Lemma 8.3]{HKM2} implies the existence of a nonnegative decreasing sequence $\{v_{i}\}\subset H_{\alpha}^{+}(U(x_{0}))$ satisfying
$$
\begin{cases}
v_{i}(x)=u(x)\quad\hbox{for all} \quad x\in E\cap V(x_{0});\\
\widehat{v}(x)=\liminf_{y\longrightarrow  x}v(y)=\widehat{\mathcal{B}}^{u}_{E\cap V(x_{0})}(x)\quad\hbox{for}\quad v=\lim_{i}v_{i}\quad\hbox{and}\quad x\in  U(x_{0}).
\end{cases}
$$
Denote by $\widetilde{E}=\left\{x\in V(x_{0}):v(x)=u(x)\right\}.$ We may assume that $E\subset \widetilde{E}$ is a Borel set and obtain
$$
\widehat{\mathcal{B}}^{u}_{\widetilde{E}}(x_{0})\leq \widehat{v}(x_{0})=\widehat{\mathcal{B}}^{u}_{E\cap V(x_{0})}(x_{0})<u(x_{0}).
$$
Moreover, (ii) of  Lemma \ref{lemma 2.2}  allows us to assume $x_{0}\in E$. Next, we prove (i) by contradiction.

Assume that $E$ is not $(\alpha, p)-$thin at $x_{0}$. Then Lemma \ref{lemma 2.1} gives the existence of a compact set $K\subset E\cap V(x_{0})$ and $K$ is not $(\alpha, p)-$thin at $x_{0}$. Let
$$
\begin{cases}
\{\varphi_{i}\}\subset C_{c}^{\infty}(U(x_{0}))\quad\hbox{be increasing with}\quad u=\lim_{i}\varphi_{i}\quad\hbox{in}\quad K;\\
h_{i}\in H_{\alpha}(U(x_{0})\backslash K)\quad\hbox{be unique with}\quad h_{i}-\varphi_{i}\in  W^{\alpha,p}(U(x_{0})\backslash K);\\
w\in H_{\alpha}^{+}(U(x_{0}))\quad\hbox{be nonnegative with}\quad w\geq u\quad\hbox{in}\quad K;\\
\widetilde{w}=\min\{w, \sup h_{i}\}, \quad w_{0}=\min\{\widetilde{w}+\varepsilon-h_{i},0\}\quad\hbox{for all}\quad\varepsilon>0.
\end{cases}
$$
Then $u\geq h_{i}$ in $U(x_{0})\backslash K$, $\widetilde{w}\in W^{\alpha,p}(V(x_{0}))$ by \cite[Theorem 1]{KKP4} and $w_{0}\in W_{0}^{\alpha,p}(U(x_{0})\backslash K)$. We conclude from the comparision principle \cite[Lemma 2.8]{Sh} that $\widetilde{w}\geq h_{i}$ in $U(x_{0})\backslash K. $ Hence  we get
$$
u\geq \mathcal{B}^{u}_{K}=\mathcal{B}^{u\chi_{K}}(U(x_{0}))\geq h_{i}\quad\hbox{in}\quad U(x_{0})\backslash K,
$$
and finally
$$
\widehat{\mathcal{B}}^{u}_{K}(x_{0})=\sup_{r>0}\inf _{B(x_{0},r)}\mathcal{B}^{u}_{K}=\min\left\{\liminf_{y\longrightarrow  x_{0}, y\in U(x_{0})\backslash K}\mathcal{B}_{K}^{u}(y), u(x_{0})\right\}\geq \min\left\{\lim_{y\longrightarrow  x_{0}, y\in U(x_{0})\backslash K}h_{i}(y), u(x_{0})\right\}=\varphi_{i}(x_{0})$$
since $K$ is not $(\alpha, p)-$thin at $x_{0}$ (see for example \cite{KKP4}). Thus, we get
$$u(x_{0})\leq \widehat{\mathcal{B}}^{u}_{K}(x_{0})\leq \widehat{\mathcal{B}}^{u}_{E\cap V(x_{0})}(x_{0})<u(x_{0})$$
which is  a contradiction. Therefore (i) is proved.

\begin{remark}
 (iii) of Theorem \ref{Theorem 1.1} shows that $x_{0}$ is not $(\alpha, p)-$regular. In fact, if $u$ is as that in (iii), $x_{0}\in V(x)\Subset U(x_{0})$. Let
$$
\begin{cases}
\{\varphi_{i}\}\subset C_{c}^{\infty}(U(x_{0}))\,\,\hbox{be increasing with}\,\, u=\lim_{i}\varphi_{i}\quad\hbox{in}\quad\partial(V(x_{0})\cap \Omega);\\
h_{i}\in H_{\alpha}(U(x_{0})\cap \Omega)\,\,\hbox{be unique with}\,\, h_{i}-\varphi_{i}\in  W^{\alpha,p}_{0}(U(x_{0})\cap\Omega).
\end{cases}
$$
Then by the proof of (iii) $\Longrightarrow$ (i), $h_{i}\leq \widehat{\mathcal{B}}^{u}_{E\cap V(x_{0})}$ in $U(x_{0})\cap \Omega$. If $x_{0}$ is $(\alpha, p)-$regular, then
$$
\varphi_{i}(x_{0})=\lim_{x\longrightarrow  x_{0}}h_{i}(x)\leq \widehat{\mathcal{B}}^{u}_{E\cap V(x_{0})}(x_{0}).$$
By letting $i\longrightarrow \infty,$  we get $u(x_{0})\leq \widehat{\mathcal{B}}^{u}_{E\cap V(x_{0})}(x_{0}),$  which is   a contradiction.
\end{remark}

\subsection{Proof of Theorem \ref{Theorem 1.3}:   characterizations of Perron  regularity  for bounded domains}\label{s3}
We prove Theorem \ref{Theorem 1.3} by showing (i)$\Longrightarrow$ (ii)$\Longrightarrow$ (i)$\Longrightarrow$ (iii)$\Longrightarrow$ (iv)$\Longrightarrow$ (v)$\Longrightarrow$ (i). The equivalence of (i) and (ii) was given in  \cite[Theorem 26]{LL2}.

(i)$\Longrightarrow $ (iii). Let
$$
\begin{cases}
\{u_{i}\}\subset C(\mathbb{R}^{n})\,\,\hbox{be increasing},\quad i=1,2,\cdots;\\
\lim_{i}u_{i}=u\,\,\hbox{in}\,\, U\,\,\hbox{and}\quad u_{i}=0\quad\hbox{in}\quad\partial V;\\
v_{i}=\begin{cases}\overline{H}_{u_{i}}^{\alpha}(V\backslash (\overline{U}\backslash \Omega))\,\,\hbox{in}\,\, V\backslash (\overline{U}\backslash \Omega);\\
u_{i}\quad\hbox{in}\quad\overline{U}\backslash \Omega.
\end{cases}
\end{cases}
$$
It follows from (ii) of Lemma \ref{lemma 3.3}  that
$$
u_{i}(x_{0})=\lim_{x\longrightarrow  x_{0}}v_{i}(x)\leq \liminf_{x\longrightarrow  x_{0}}\mathcal{B}^{u}_{\overline{U}\backslash \Omega}(V)(x)=\widehat{\mathcal{B}}^{u}_{\overline{U}\backslash \Omega}(V)(x_{0}),
$$
and hence $u(x_{0})\leq \widehat{\mathcal{B}}^{u}_{\overline{U}\backslash \Omega}(V)(x_{0})$, which implies (iii) since it is not hard to get $u(x_{0})\geq \widehat{\mathcal{B}}^{u}_{\overline{U}\backslash \Omega}(V)(x_{0}). $

(iii)$\Longrightarrow $(iv) is straightforward.

(iv)$\Longrightarrow  $(v) follows from \cite[Proposition 3.2]{SX3}.

(v)$\Longrightarrow  $(i). This can be done by showing (v)$\Longrightarrow $ (iv)$\Longrightarrow$ (i). Assume that $x_{0}$ is $(\alpha,p)-$regular and that $B$ is a ball containing $x_{0}$. Then \cite[Theorem 1.2]{SX3} implies
$$
\widehat{\mathcal{B}}^{1}_{\overline{B}\backslash \Omega}(2B)(x_{0})=\lim_{x\longrightarrow  x_{0}, x\in \Omega}\widehat{\mathcal{B}}^{1}_{\overline{B}\backslash \Omega}(2B)(x)=1,
$$
which is (iv). (iv)$\Longrightarrow $ (i) is immediately from Lemma \ref{lemma 3.2}. This finishes the proof of Theorem \ref{Theorem 1.3}.

{\bf The regularity is a local property}.
By the equivalence (i)$\Longleftrightarrow$(ii) of Theorem \ref{Theorem 1.3}, we get the following fact that regularity is  a local property.
\begin{corollary}\label{lemma 3.3}
Let $O, \Omega \subset \mathbb{R}^{n}$ be open sets, $x_{0}\in \partial\Omega \cap \partial O$ and $U(x_{0})$ be a neighborhood of $x_{0}$.
\begin{itemize}
\item[\rm (i)]  If $O\subset \Omega$ such that  $U(x_{0})\cap O=U(x_{0})\cap \Omega$, then
there exists an $(\alpha, p)-$barrier relative to $\Omega$  at $x_{0}$ if and only if there exists an $(\alpha, p)-$barrier relative to $O$ at $x_{0}.$

\item[\rm (ii)] If $U(x_{0})\cap O=U(x_{0})\cap \Omega$, then
$
x_{0}$ is Perron regular relative to $\Omega$ if and only if $x_{0}$ is Perron regular relative to $ O.$

\item[\rm (iii)] Assume that $O\subset\Omega\subset \mathbb{R}^{n}$. Then $x_{0}$ is Perron regular relative to $O$ if $x_{0}$ is Perron regular relative to $\Omega$.
\end{itemize}
\end{corollary}

\begin{proof}
The equivalence of (i) and (ii) in  Theorem \ref{Theorem 1.3} shows that the existence of an $(\alpha, p)$-barrier is a local property, i.e. (i), (ii) is a consequence of (i)$\Longleftrightarrow$(ii) of Theorem \ref{Theorem 1.3} and (i), while (iii) is a byproduct of (ii).
\end{proof}

\subsection{Proof of Theorem \ref{Theorem 1.4}:  characterizations of Perron regularity for unbounded domains}\label{s4}

To prove Theorem \ref{Theorem 1.4}, we need  the following lemma concerning the weak compactness in $\dot{W}^{\alpha,p}_{0}(\Omega)$.
\begin{lemma}\label{lemma 4.1}
Let $\{u_{i}\}\subset \dot{W}^{\alpha,p}_{0}(\Omega)$ be a  bounded squence  in $\dot{W}^{\alpha,p}(\Omega)$ satisfying  $\lim_{i}u_{i}=u$ and
$
u_{i}$ being bounded in $\dot{W}^{\alpha,p}(O)$ for each open set $ O\Subset \Omega.$
Then
$
u\in \dot{W}^{\alpha,p}_{0}(\Omega)$ and
$u_{i}\longrightarrow  u\,\,\hbox{weakly in}\,\,\dot{W}^{\alpha,p}(\Omega).
$
\end{lemma}

\begin{proof}
From \cite[Lemma 2.2]{Sh}, it follows that $u\in \dot{W}^{\alpha,p}_{loc}(\Omega)$ and $u_{i}\longrightarrow  u$ weakly in $\dot{W}^{\alpha,p}(O)$ for $O\Subset \Omega$. The $\dot{W}^{\alpha,p}(\Omega)$ boundedness of $\{u_{i}\}$ gives the existence of a weakly convergent subsequence $\{u_{i_{j}}\}$ of $\{u_{i}\}$. Hence, the weak convergence of  $u_{i_{j}}\longrightarrow  u$   implies $u\in \dot{W}^{\alpha,p}(\Omega).$ Furthermore, we have $u_{i}\longrightarrow  u$ weakly in $\dot{W}^{\alpha,p}(\Omega)$ since the independence of the subsequence for the weak limit. We conclude from the Mazur lemma \cite[Lemma 1.29]{HKM2} that
$u\in \dot{W}^{\alpha,p}_{0}(\Omega)$. In fact, by choosing a sequence $\{v_{i}\}\subset \dot{W}^{\alpha,p}_{0}(\Omega)$ of convex combinations of $u_{i}$ with $v_{i}\longrightarrow  u$ in $\dot{W}^{\alpha,p}(\Omega)$ and
$$[v_{i}-u]_{\dot{W}^{\alpha,p}(\Omega)}<\left(\frac{\varepsilon}{2}\right)^{p}, \quad[v_{i}-\varphi]_{\dot{W}^{\alpha,p}(\Omega)}<\left(\frac{\varepsilon}{2}\right)^{p}\quad\hbox{for}\quad  \varphi\in C_{c}^{\infty}(\Omega),\quad\varepsilon>0.$$
 Thus, we achieve the desired result
$[\varphi-u]_{\dot{W}^{\alpha,p}(\Omega)}<\varepsilon$.
\end{proof}
{\bf Proof of Theorem \ref{Theorem 1.4}}.    Theorem \ref{Theorem 1.4} will be proved by showing
(i)$\Longrightarrow $(ii)$\Longrightarrow $(iii)$\Longrightarrow $(iv)$\Longrightarrow $(v)$\Longrightarrow $(i).

(i)$\Longrightarrow $ (ii) is straightforward.

(ii)$\Longrightarrow $ (iii). Assume that $\infty$ is Perron regular for $\Omega=\overline{B}^{c}$ with $B$ being a ball. Then $\Omega$ is regular since $\overline{B}$ is $(\alpha,p)-$thick at each of its points. This  implies the existence of a function $u\in H_{\alpha}(\Omega)$ such that $u\in C(\overline{\Omega})$, $u=1$ on $\partial B$ and $u=0$ at $\infty$. Consequently, $u\in H_{\alpha}^{+}(\mathbb{R}^{n})$ is the desired nonconstant bounded function with $0\leq u\leq1$ by extending $u$ to $B$ as $u=1$ in $B$.

(iii)$\Longrightarrow $ (iv). Let $u$ be the function in (iii). Without loss of generality, we assume $\inf u=0$. Choosing an open ball $B\subset \mathbb{R}^{n}$, we observe that $m=\min _{\overline{B}}u>0.$ We complete the proof by a contradiction.

Suppose that $C_{\alpha,p}(\overline{B}, \mathbb{R}^{n})=0.$ Then,  $\mathbb{R}^{n}$ can be exhausted by increasing concentric balls $\{B_{i}\}$ with $B\Subset B_{1}\Subset B_{2}\Subset\cdots $ and $C_{\alpha,p}(\overline{B}, B_{i})\longrightarrow  0$.
Let $u_{i}=\widehat{\mathcal{B}}^{1}_{\overline{B}}(B_{i})$. Then $u_{i}\in H_{\alpha}(B_{i}\backslash \overline{B})$ and $u_{i}\longrightarrow  v\in H_{\alpha}(\mathbb{R}^{n}\backslash \overline{B})$ by \cite[Theorem 1.2]{SX3} and \cite[Theorem 15]{KKP4}. It follows from the fact $v\leq u/m$ that $v$ is not constant. In fact, $\lim_{x\longrightarrow  y}h(x)=1$ for $y\in \partial B$. At the same time, the quasiminimizing property \cite[P 1056]{Sh} gives $[u_{i}]_{\dot{W}^{\alpha,p}(B_{i})}\leq C_{\alpha,p}(\overline{B}, B_{i})$. Namely, $\{u_{i}\}$ is bounded in $\dot{W}^{\alpha,p}(\mathbb{R}^{n})$. We conclude from Lemma \ref{lemma 4.1} that $v\in \dot{W}^{\alpha,p}(\overline{B}^{c})$ and $u_{i}\longrightarrow  v$ weakly in $\dot{W}^{\alpha,p}(\overline{B}^{c}).$ Hence  it follows from \cite[Remark 5.25]{HKM2} that
$$
[v]_{\dot{W}^{\alpha,p}(\overline{B}^{c})}\leq \liminf_{i\longrightarrow \infty}[u_{i}]_{\dot{W}^{\alpha,p}(\overline{B}^{c})}\leq \lim_{i\longrightarrow \infty}C_{\alpha,p}(\overline{B}, B_{i})=0.
$$
Finally, we can derive that $v$ is constant in $\overline{B}^{c}$, which is a contradiction and thus (iv) follows.

(iv)$\Longrightarrow $(v) is straightforward.

(v)$\Longrightarrow $(i). Let
$$
\begin{cases}
C_{\alpha,p}(\overline{B}, \mathbb{R}^{n})=\delta>0\quad\hbox{for a ball}\quad B;\\
u_{i}=\widehat{\mathcal{B}}^{1}_{\overline{B}}(B_{i})\quad\hbox{with}\quad B_{i}=(i+1)B,\quad i=1,2,\cdots.
\end{cases}
$$
Then $u_{i}\longrightarrow  u\in H_{\alpha}(\overline{B}^{c}).
$
We claim first that $u$ is not constant. In fact, by the same analysis as in the proof of (iii)$\Longrightarrow  $(iv), we have $u_{i}\longrightarrow  u$ weakly in $\dot{W}^{\alpha,p}(\overline{B}^{c})$. The Mazur lemma \cite{HKM2} allows the existence of a convex combination of $u_{i}$, $$\overline{u_{i}}=\sum_{j=1}^{i}\lambda_{j,i}u_{j},\quad \sum_{j=1}^{i}\lambda_{j,i}=1\quad \hbox{ with}\quad  \lambda_{j,i}\geq 0$$
and $\overline{u_{i}}\longrightarrow  u$ in $\dot{W}^{\alpha,p}(\overline{B}^{c})$. If $u$ is a constant, then there exists a $i_{k}$ with $[\overline{u_{i_{k}}}]_{\dot{W}^{\alpha,p}(\overline{B}^{c})}<\delta$. However, $C_{\alpha,p}(\overline{B}, \mathbb{R}^{n})\leq C_{\alpha,p}(\overline{B}, B_{i_{k}})<\delta.$ This  forces a contradiction.

Next, we claim that $\infty$ has a strictly positive $(\alpha,p)-$barrier relative to $\mathbb{R}^{n}$. Indeed, by \cite[Lemma 2.5]{SX3}, we infer that $\lim_{x\longrightarrow  \infty}u(x)=\inf u<1.$ Moreover, $u\in H_{\alpha}^{+}(\mathbb{R}^{n})$ if $u=1$ on $\overline{B}$ since $\lim_{x\longrightarrow  y\in\partial B}u(x)=1$. Accordingly, $v=u-\inf_{u}$ is the desired $(\alpha,p)-$barrier, which implies that $v$ is an $(\alpha,p)-$barrier at $\infty$ relative to each unbounded open set. This  yields (i) by Theorem \ref{Theorem 1.3}. Then the proof of Theorem \ref{Theorem 1.4} is finished.

 If $\Omega$ is bounded in Theorem \ref{Theorem 1.4}, we have the following Corollary.
\begin{corollary}\label{Corollary 1.1}
Assume that $\infty$ is a Perron irregular boundary point of an unbounded open set $E$ and $\Omega\subset \mathbb{R}^{n}$ is a bounded open set. Then the following statements are equivalent:
\begin{itemize}	
\item[\rm (i)] $C_{\alpha, p}(\Omega^{c})>0$;

\item[\rm (ii)] There exists a bounded $u\in H_{\alpha}^{+}(\Omega)$ satisfying that $u$ is nonconstant in each component of $\Omega$;

\item[\rm (iii)] $C_{\alpha, p}(B, \Omega)>0$ for each ball $B\Subset \Omega$;

\item[\rm (iv)] $C_{\alpha, p}(B, \Omega)>0$ for some ball $B\Subset \Omega$.
\end{itemize}
\end{corollary}
 \begin{proof}
 (i)$\Longrightarrow $(ii)$\Longrightarrow $(iii)$\Longrightarrow $ (iv) go essentially as that in the proof of Theorem \ref{Theorem 1.4} due to the fact that there exists at least one Perron regular point of $\Omega\backslash \overline{B}$ on $\partial \Omega$ if $B\Subset \Omega$ is a ball, and that the balls $\{B_{i}\}$ in the  argument can be replaced by polyhedra.

The only difference lies in showing (iv)$\Longrightarrow $ (i), which can be done by contradiction. Since $\infty$ is a Perron irregular boundary point of $\Omega$, Theorem \ref{Theorem 1.4} gives $C_{\alpha,p}(\overline{B}, \mathbb{R}^{n})=0.$ Then,  for fixed $\varepsilon>0$, we can pick a function $\varphi\in Y^{\alpha,p}_{0}(\overline{B}, \mathbb{R}^{n})$ such that $[\varphi]_{\dot{W}^{\alpha,p}(\mathbb{R}^{n})}<\varepsilon$. Assume  $C_{\alpha,p}(\Omega^{c})=0$ and let $K=\supp \varphi\cap \Omega^{c}.$ Then we get $C_{\alpha,p}(K)=0$. We thus pick a function $\psi\in C_{c}^{\infty}(\mathbb{R}^{n})$ with
$$
0\leq \psi \leq 1;\quad
\psi=1\quad\hbox{in a neighborhood of}\quad K;\quad
\psi=0\quad\hbox{on}\quad\overline{B};\quad
[\psi]_{\dot{W}^{\alpha,p}(\mathbb{R}^{n})}<\varepsilon.
$$
Consequently, we get $\eta=(1-\psi)\varphi\in Y^{\alpha,p}_{0}(\overline{B}, \Omega).$ Hence, we derive
$$
C_{\alpha,p}(\overline{B},\Omega)\leq [\eta]_{\dot{W}^{\alpha,p}(\Omega)}\leq 2^{p-1}\left([\varphi]_{\dot{W}^{\alpha,p}(\mathbb{R}^{n})}+[\psi]_{\dot{W}^{\alpha,p}(\mathbb{R}^{n})}\right)\leq 2^{p}\varepsilon,
$$
which is a contradiction of (iv) and then (i) follows.
\end{proof}

\subsection{Proof of Theorem \ref{Theorem 1.2}: the fractional Wiener test} \label{s5}

In this section, we show the fractional Wiener test by using Theorem \ref{Theorem 1.1}.

The second equivalence of Theorem \ref{Theorem 1.2} follows immediately by the definition of $(\alpha, p)$-thickness. We only need to show that the $(\alpha, p)$-regular boundary points can be characterized by the fractional Wiener integral. The regularity of $x_{0}$ can be achieved by the divergence of the fractional Wiener integral according to \cite[Theorem 1.1]{SX3}. It suffices to show that the regularity of $x_{0}\in \partial\Omega$ implies the divergence of $\mathcal{W}^{\alpha}_{p}(\Omega^{c}, x_{0})$. To do this,  we need to prove the following lemmas which are critical  not only for the proof of Theorem \ref{Theorem 1.2} but also for  other results in Section \ref{s8}.

\subsubsection{Important Lemmas} It is well known that functions in $H_{\alpha}^{+}(\Omega)$ are not local since it requires testing in all open sets $O\Subset \Omega$. However, we can show the local nature of $(\alpha, p)-$superharmonic functions by \cite[Proposition 2.23, Theorem 3.6, Corollarys 3.7\ \&\ 3.9]{Sh}. Indeed, we have the following result.
\begin{lemma}\label{lemma 3.1}
	Let $\Omega$ be an open set. Then
\begin{itemize}
\item[\rm (i)]  $u\in H_{\alpha}^{+}(\Omega)\Longleftrightarrow u|_{U}\in H_{\alpha}^{+}(\Omega)$ with $U\subset \Omega$ be the neighborhood of $x\in \Omega$.

\item[\rm (ii)] Let
$$
\begin{cases}
O_{i},\quad i=1,2\quad\hbox{be open sets with}\quad\Omega=O_{1}\cup O_{2};\\
u\in H_{\alpha}^{+}(O_{1}),\quad v\in H_{\alpha}^{+}(O_{2})\quad\hbox{satisfying}\quad u\leq v\quad\hbox{in}\quad O_{1}\cap O_{2}\\
\quad \hbox{and}\quad w=\begin{cases}v\quad\hbox{on}\quad\Omega\backslash O_{1};\\
u\quad\hbox{in}\quad O_{1}
\end{cases}
is\quad l.s.c.
\end{cases}
$$
Then $w\in H_{\alpha}^{+}(\Omega)$.
\end{itemize}
\end{lemma}
\begin{proof}
Since we may assume that $u$ is bounded,
it remains to consider the corresponding property of supersolutions according to  \cite[Proposition 2.16, Theorems 3.1 \& 3.6]{Sh}.   Therefore, (i) follows by applying a partition of unity.

Using (i), we can prove (ii), which can be seen as another version of the classical pasting lemma \cite[Proposition 2.21]{Sh}. In fact, it follows from \cite[Proposition 2.21]{Sh} that $w\in H_{\alpha}^{+}(O_{2})$ since $\min\{u,v\}=u$ in $O_{1}\cap O_{2}$. Hence, $w\in H_{\alpha}^{+}(\Omega)$ by (i) noticing that $O_{1}\cup O_{2}=\Omega.$
\end{proof}

\begin{lemma}\label{lemma 3.2}
Let $O$ be an open set and $x_{0}\in B\cap (\partial O \backslash \{\infty\})$ for a ball $B$ with rational center and radius. Then

\begin{itemize}
\item[\rm (i)] $x_{0}$ is Perron regular if  $\widehat{\mathcal{B}}^{1}_{\overline{B}\backslash O}(2B)(x_{0})=1$.

\item[\rm (ii)] $C_{\alpha,p}(S,O)=0$, where $S=\{x:x\in \partial O \,\,\hbox{is finite and}\,\,(\alpha,p)-\hbox{irregualr}\}$.

\item[\rm (iii)] There exists a finite point $x_{0}\in \partial O$ which is  Perron regular relative to $E\subset O$ if $C_{\alpha,p}(O^{c})>0$.

\item[\rm (iv)] Let $u\in H_{\alpha}^{+}(U)$  and $C_{\alpha,p}(O^{c})>0$ with $U$ be a neighborhood of $\overline{B}\subset O$. Then there exists a function $\widetilde{u}\in H_{\alpha}^{+}(O)$ with $\widetilde{u}=u$ in $B$ and $\widetilde{u}$ is bounded below.
\end{itemize}
\end{lemma}
\begin{proof}
(i) is a slight modification of \cite[Proposition 3.2]{SX3} with  $v$ and $h$  replaced by $f\in C(\partial O)$, $\overline{H}_{f}^{\alpha}$.

(i)$\Longrightarrow$(ii) can be derived in a similar way as that of (ii) $\Longrightarrow$ (iii) in  \cite[Proposition 3.2]{SX3}.

To prove (iii), we first claim that $C_{\alpha,p}(\partial E)>0$ since $C_{\alpha,p}(E^{c})=0$. In fact, if $C_{\alpha,p}(\partial E)=0$, then $E^{c}$ must have an interior point. This  implies that $(\partial E)^{c}$ has at least two components, which is a contradiction since $\partial E$ does not separate $\mathbb{R}^{n}$(see for example \cite[Lemma 2.46]{HKM2}, \cite[Proposition 4.4]{SX1} and  \cite[Lemma 3.3]{SZ}). Consequently, $C_{\alpha,p}(\partial E)>0$ and (iii)
follows from (ii).

(iv). Without loss of generality, we may assume that $O$ is connected since it is sufficient to prove that $u$ has a desired extension to the component of $O$ which contains $B$. Pick a ball $B_{0}\Subset O$ satisfying
$$
u\in H_{\alpha}^{+}(B_{0});\quad
B\Subset B_{0};\quad
u>0\quad\hbox{in}\quad B_{0}
$$
and replace $u$ instead of $v=\widehat{\mathcal{B}}^{u}_{B}(B_{0}).$ It follows from  \cite[Lemma 2.7 \& Proposition 3.2]{SX1} that
$$
v=u\quad\hbox{in}\quad B\quad\&\quad\lim_{x\longrightarrow  y}v(x)=0\quad\hbox{for  all}\quad y\in \partial B_{0}.
$$
Set $m=\min_{\overline{B}}v, $ $ K=\left\{x\in B_{0}: v(x)\geq m\right\}\supset \overline{B},
$
$$
w=\begin{cases}\overline{H}_{f}^{\alpha}\quad\hbox{in}\quad O\backslash K;\\
m\quad\hbox{in}\quad K;
\end{cases}\quad \hbox{and}\quad
f=\begin{cases}0\quad\hbox{on}\quad\partial O;\\
m\quad\hbox{in}\quad K.
\end{cases}\\
$$
We proceed the proof by showing
\begin{equation}\label{3.2}
\lim_{x\longrightarrow  y}w(x)=m\quad \forall y\in K.
\end{equation}
Choose $\widetilde{w}\in \mathcal{U}_{f}^{\alpha}$, then \cite[Proposition 2.20]{Sh} implies that $\widetilde{w}$ is nonnegative and $\widetilde{w}\geq v$ in $B_{0}\backslash K$ since $v\in H_{\alpha}(B_{0}\backslash \overline{B})$. Consequently, $\lim_{x\longrightarrow  y}w(x)\geq m$  $\forall y\in \partial K$ and $(\ref{3.2})$ is obtained. Furthermore, \cite[Proposition 2.21]{Sh} implies $w\in H_{\alpha}^{+}(O)\cap H_{\alpha}(O\backslash K).$ Then (iii) yields the existence of a point $y\in \partial O$ with $\lim_{x\longrightarrow  y}w(x)=0$, which shows  $w<m$ on $\partial B_{0}$ by \cite[Lemma 2.5]{SX3}. Thus, we have
\begin{equation}\label{3.3}
m/((m-\max_{\partial B_{0}}w)(w-m))\leq -m=v-m\quad\hbox{on}\quad\partial B_{0}.
\end{equation}
Similarly, a further application of \cite[Proposition 2.20]{Sh} implies the truth of   $(\ref{3.3})$  in $B_{0}\backslash K$. Finally, Lemma \ref{lemma 3.1} implies that the function
$$
\widetilde{u}=\begin{cases}
v\quad\hbox{in}\quad K;\\
(m-\max_{\partial B_{0}}w)(w-m)+m\quad\hbox{in}\quad O\backslash K
\end{cases}
$$
is the desired function.
\end{proof}

To prove Theorem \ref{Theorem 1.2}, we  need the continuity of $(\alpha, p)$-balayage which can be formulated as follows.
\begin{lemma}\label{lemma 3.4}
Assume that $f:\overline{\mathbb{R}^{n}}\longrightarrow  \mathbb{R}^{n}$ is continuous and $u=\widehat{\mathcal{B}}^{f}_{\Omega}$. Then
\begin{itemize}
\item[\rm (i)]  $u$ is continuous, $u\geq f$ and $\lim_{x\longrightarrow  x_{0}}u(x)=f(x_{0})$ if $x_{0}\in \partial\Omega\backslash\{\infty\}$ is Perron regular.

\item[\rm (ii)] $u\in \dot{W}^{\alpha,p}(\Omega)$, $u-f\in \dot{W}_{0}^{\alpha,p}(\Omega)$ and $\langle(-\Delta_{p})^{\alpha}u, v\rangle\geq 0$ if
$f\in \dot{W}^{\alpha,p}_{0}(\Omega)$ and $\varphi\in \dot{W}_{0}^{\alpha,p}(\Omega)$ satisfy $\varphi\geq f-u$.

\item[\rm (iii)] $u\in W^{\alpha,p}(\Omega)$, $u-f\in W_{0}^{\alpha,p}(\Omega)$ if
$f\in W^{\alpha,p}(\Omega)$ and $\Omega$ is bounded.

\item[\rm (iv)] $\lim_{x\longrightarrow  \infty}u(x)=f(\infty)$ if $\Omega$ is unbounded and $\infty$ is a Perron regular boundary point of each unbounded open set $O\subset \Omega$.
 \end{itemize}
\end{lemma}
\begin{proof}
It follows from \cite[Proposition 3.3]{SX3} that $u\in C(\Omega)$ is bounded and $u\geq f$. Suppose that $x_{0}\in \partial\Omega\backslash \{\infty\}$ is Perron regular. We proceed to show that
\begin{equation}\label{3.1}
\lim_{x\longrightarrow  x_{0}}u(x)=f(x_{0}).
\end{equation}
By letting
$$
\begin{cases}
v\,\,\hbox{be an}\,\,(\alpha,p)-\hbox{barrier at}\,\, x_{0};\\
B\,\,\hbox{be a ball centered at}\,\, x_{0}\,\,\hbox{with}\,\,\partial B\cap\Omega \neq \phi;\\
|f(x)-f(x_{0})|<\varepsilon\,\,\hbox{for all}\,\, x\in B\,\,\hbox{and fixed}\,\,\varepsilon>0,
\end{cases}
$$
we obtain from \cite[Proposition 2.21]{Sh} the existence of a $\xi>0$ and
$$
w(x)=\begin{cases}
\min\{\xi v(x)+f(x_{0})+\varepsilon, \sup |f|\}\quad\hbox{for}\quad x\in B\cap \Omega;\\
\sup |f|\quad\hbox{for}\quad x\in  \Omega\backslash B,
\end{cases}
$$
which is  l.s.c   satisfying $w\in H_{\alpha}^{+}(\Omega)$ and $w\geq f$ in $\Omega$. Hence,  we get
$$\limsup_{x\longrightarrow  x_{0}}u(x)\leq \lim_{x\longrightarrow  x_{0}}w(x)\leq f(x_{0})+\varepsilon,$$ which implies $(\ref{3.1})$ since $u\geq f$.

To show (ii), let
$$
\begin{cases}
f\in \dot{W}^{\alpha,p}(\Omega);\\
\{O_{i}\}\,\,\hbox{be increasing regular open sets with}\,\, O_{1}\Subset O_{2}\Subset\cdots \Subset \Omega\,\,\hbox{and}\,\, \cup_{i}O_{i}=\Omega;\\
u_{i}\,\,\hbox{be the solution to the obstacle problem}\,\,\Phi_{f,f}(O_{i}).
\end{cases}
$$
It follows from  \cite[Lemma 2.6]{Sh}  that $\{u_{i}\}$ is increasing and bounded above by $u$. Therefore, $$\lim_{i}u_{i}=\widetilde{u}\in H_{\alpha}^{+}(\Omega)\quad \hbox{with}\,\, u\geq \widetilde{u}\geq f.$$ Consequently, we derive  $u=\widetilde{u}$ since $u\leq \widetilde{u}$ by the definition of $(\alpha,p)-$balayage.

By  letting $u_{i}=f$ on $\Omega\backslash O_{i}$, we see that $u_{i}-f\in \dot{W}^{\alpha,p}_{0}(\Omega)$. On the other hand, the quasi minimizing property (see for example \cite[P1056]{Sh}) implies that
$$
[u_{i}]_{\dot{W}^{\alpha,p}(O_{i})}\leq C[f]_{\dot{W}^{\alpha,p}(O_{i})}\leq C[f]_{\dot{W}^{\alpha,p}(\Omega)}.
$$
Accordingly, we deduce that $(u_{i}-f)\longrightarrow  (u-f)$ weakly in $\dot{W}^{\alpha,p}(\Omega)$, and hence $(u-f)\in \dot{W}^{\alpha,p}_{0}(\Omega)$ by \cite[Lemma 2.2]{Sh}. Let $\{\varphi_{i}\}\subset C_{c}^{\infty}(\Omega)$ such that $\varphi_{i}\longrightarrow  \varphi$ in $\dot{W}^{\alpha,p}(\Omega)$. Then the function $\psi_{i}=\max\{\varphi_{i}, f-u\}\in \dot{W}^{\alpha,p}_{0}(\Omega),$ and $\psi_{i}\longrightarrow  \varphi$ in $\dot{W}^{\alpha,p}(\Omega)$ by \cite[Lemma 2.1]{SX3}. For fixed $i$, we choose $j_{i}$ with $\supp \psi_{i}\subset O_{j_{i}}$. It follows from \cite[Proposition 2.15]{Sh} that $u$ is a solution to $\Phi_{f,u}(O_{j_{0}})$. Thus,
$$
\langle(-\Delta_{p})^{\alpha}u, \psi_{i}\rangle=\langle(-\Delta_{p})^{\alpha}u, \psi_{i}\rangle|_{O_{j_{i}}}\geq 0.
$$
Consequently, we obtain
$$
\langle(-\Delta_{p})^{\alpha}u, \varphi\rangle=\lim_{i\longrightarrow  \infty}\langle(-\Delta_{p})^{\alpha}u, \psi_{i}\rangle\geq 0,
$$
which complete the proof of (ii).

(iii) and (iv) is a slight modification of (ii) and (i), respectively.
\end{proof}

In fact, the continuity assumption in Lemma \ref{lemma 3.4} can be relaxed to an extention in the following sense, which  can be established  in a similar way as that of (i) and (ii) of Lemma \ref{lemma 3.4}.

\begin{corollary}\label{corollary 3.5}
Let $f$ be bounded above, $u=\widehat{\mathcal{B}}^{f}_{\Omega}$ and $x_{0}\in \partial\Omega\backslash \{\infty\}$ be Perron regular. Then
\begin{itemize}
\item[\rm (i)]  $\lim_{x\longrightarrow  x_{0}}u(x)=f(x_{0})$ whenever $f$ is continuous at $x_{0}$ and $u\geq f$ in $B(x_{0},r)\cap\Omega$ for some $r>0$.

\item[\rm (ii)] $u\in \dot{W}^{\alpha,p}(\Omega)$, $u-f\in \dot{W}_{0}^{\alpha,p}(\Omega)$ and $\langle(-\Delta_{p})^{\alpha}u, v\rangle\geq 0$ if
$f\in \dot{W}_{0}^{\alpha,p}(\Omega)$ is bounded above and the solution $v$ to $\Phi_{f,f}(O)$ lies above $f$ in all open set $O\Subset \Omega$.
\end{itemize}
\end{corollary}

Utilizing Lemma \ref{lemma 3.4}, one can get a solution to $(\ref{(1.2)})$ with Sobolev boundary values, which helps us to complete the proof of Theorem \ref{Theorem 1.2}.
\begin{lemma}\label{lemma 3.6} The following two statements hold.
\begin{itemize}
\item[\rm (i)] Assume that $C_{\alpha,p}(\Omega^{c})>0$ and that $f\in C(\overline{\mathbb{R}^{n}})\cap \dot{W}^{\alpha,p}(\Omega)$. Then $f$ is resolutive and $\overline{H}^{\alpha}_{f}-f\in \dot{W}_{0}^{\alpha,p}(\Omega)$ if $\Omega$ is an unbounded open set.

\item[\rm (ii)] Assume that $f\in C(\overline{\mathbb{R}^{n}})\cap \dot{W}^{\alpha,p}(\Omega)$. Then $\overline{H}^{\alpha}_{f}-f\in W_{0}^{\alpha,p}(\Omega)$ if $\Omega$ is bounded. Moreover,  $\overline{H}^{\alpha}_{f}$ is the unique $(\alpha,p)-$harmonic function with Sobolev boundary values $f$.
\end{itemize}
\end{lemma}
\begin{proof}
Let $u=\widehat{\mathcal{B}}^{f}_{\Omega}$ and $O_{i}$ be as in the proof of Lemma \ref{lemma 3.4} (ii). Denote by $u_{i}$  the Poisson modification $P(u,O_{i})$ of $u$ in $O_{i}$ (see \cite[P1069]{Sh}). It is easy to see that
$$u\geq u_{1}\geq u_{2}\geq \cdots\geq \overline{H}_{f}^{\alpha}\quad\hbox{and}\quad  \overline{H}_{f}^{\alpha}\leq \lim_{i}u_{i}=u^{*}\in H_{\alpha}(\Omega)$$
according to \cite[Proposition 2.22 \& Remark 2.1]{Sh}. Furthermore, we get $u_{i}\longrightarrow  u^{*}$ weakly in $\dot{W}^{\alpha,p}(\Omega)$ and $f-u^{*}\in \dot{W}_{\alpha,p}(\Omega)$ by a slight modification of the proof of Lemma \ref{lemma 3.4} (i) and (ii).

Similarly, the $(\alpha,p)-$subharmonic function $v=-\widehat{\mathcal{B}}^{-f}_{\Omega}$ gives the existence of a function $u_{*}\in H_{\alpha}(\Omega)$ such that
$v\leq u_{*}\leq \underline{H}_{f}^{\alpha}\leq  \overline{H}_{f}^{\alpha}\leq u^{*}$ in $\Omega$ and $f-u_{*}\in \dot{W}_{0}^{\alpha,p}(\Omega)$. To complete the proof of (i), it remains to verify $u^{*}= u_{*}$. Since $u^{*}- u_{*}\in \dot{W}_{0}^{\alpha,p}(\Omega)$, we get
$
\langle(-\Delta_{p})^{\alpha}u^{*}-(-\Delta_{p})^{\alpha}u_{*}, u^{*}- u_{*}\rangle=0,
$
and hence $u^{*}= u_{*}+C$ for some constants $C$ in each component $\Omega_{0}$ of $\Omega$. Next, we claim that $C=0.$ In fact, Lemma \ref{lemma 3.2} (iii) shows that $\Omega_{0}$ has a finite Perron regular boundary point $x_{0}$ since $C_{\alpha,p}(\Omega^{c})>0$, and hence Lemma \ref{lemma 3.4} now implies, as   $x\in \Omega_{0}$ goes to $  x_{0},$
$$
\limsup u_{*}(x)\leq \limsup u^{*}(x)\leq \lim u(x)=f(x_{0})=\lim v(x)\leq \liminf u_{*}(x)\leq \liminf u^{*}(x),
$$
which is the desired result.

(ii) can be proved similarly as that of (i). We may assume that $f\in C(\overline{\mathbb{R}^{n}})$ by the Tietze extension theorem. Then $(u_{i}-f)\longrightarrow  \overline{H}_{f}^{\alpha}-f$ weakly in $W_{0}^{\alpha,p}(\Omega)$ as $\widehat{\mathcal{B}}^{f}-f\in W_{0}^{\alpha,p}(\Omega)$ by Lemma \ref{lemma 3.4} (iii) and (iv).

On the other hand, (ii) can also be proved by the comparision principle. Let $h\in H_{\alpha}(\Omega)$ such that $h-f\in W_{0}^{\alpha,p}(\Omega)$. Then $\lim_{x\longrightarrow  y}h(x)=f(y)$ for $x\in \partial\Omega$ q.e. by Lemma \ref{lemma 3.2}. Then, it follows from \cite[Proposition 2.20]{Sh}  that $h\leq u$ for $u\in \mathcal{U}_{f}^{\alpha}$. On the other hand, $h\geq v$ for $v\in \mathcal{L}_{f}^{\alpha}$, which shows that $h\leq \overline{H}_{f}^{\alpha}=\underline{H}_{f}^{\alpha}\leq h$ and finishes the proof.
\end{proof}

\subsubsection{Proof of Theorem \ref{Theorem 1.2}} Having established the previous  technical lemmas, we  are ready to prove Theorem \ref{Theorem 1.2}. Next, we show that the regularity of $x_{0}$ implies the divergence of the Wiener type integral by contradiction. Assume that $W_{\alpha,p}(\Omega^{c}, x_{0})< \infty$, then the proof follows by two cases.

Case 1.  $x_{0}$ is an isolated boundary point. In this case, it is easy to show that $x_{0}$ is irregular by \cite[Lemma 2.5]{SX3} and \cite[Theroem 3.4]{SZ}.

Case 2.  $x_{0}$ is an accumlation point of $\Omega^{c}$. Using Theorem \ref{Theorem 1.1}, we deduce that
there exist  balls  $B_{i}:=B(x_{0},r_{i}),  i=1,2,$ with $r_{1}<r_{2},$
and a function  $u\in H_{\alpha}^{+}(B_{2})$ with $ 0\leq u\leq 1,$ $u=1$ in $B_{2}\cap (\Omega^{c}\backslash \{x_{0}\})$ and $u(x_{0})\leq \frac{1}{2}.$
Pick a function $\varphi\in C^{\infty}(\mathbb{R}^{n})$ satisfying
$$
\begin{cases}
\varphi(x)\leq u(x)\quad\forall x\in \Omega^{c}\cap \overline{B}_{1}\backslash \{x_{0}\};\\
\varphi(x)=1\quad\forall x\in U(x_{0})(\hbox{a neighborhood of}\, x_{0}).
\end{cases}
$$
Then, we can derive  $\overline{H}_{\varphi}^{\alpha}\in W^{\alpha,p}(B_{1}\cap \Omega)$  from Lemma \ref{lemma 3.6} (ii) and the fact that the set of the irregular boundary points is of $(\alpha,p)$-capacity zero \cite[Proposition 3.2]{SX3}.
Thus,
we get  $\overline{H}_{\varphi}^{\alpha}(x)\leq u(x)$ for $x\in B_{r}\cap \Omega$ according to
 \cite[Corollary 2]{KKP4}.
Particularly,
$$
\liminf_{x\longrightarrow  x_{0}}\overline{H}_{\varphi}^{\alpha}(x)\leq \liminf_{x\longrightarrow  x_{0}, x\in \Omega}u(x)=u(x_{0})\leq \frac{1}{2}<1=\varphi(x_{0}).
$$
Accordingly, $x_{0}$ is not $(\alpha,p)$-regular boundary point of $B_{1}\backslash \Omega$, which shows that $x_{0}$ is not an $(\alpha,p)$-regular boundary point of $\Omega$ by the fact that the regularity is a local property(Lemma \ref{lemma 3.3}). The proof of Theorem \ref{Theorem 1.2} is completed.

\subsection{Proof of Theorem \ref{Theorem 5.1}: the generalized fractional Wiener criterion}\label{sec2.5}

Assume that $u\in W^{\alpha,p}(\Omega)$ is a local solution to $\Phi^{x_0}_{f,u_{0}}.$ Then $u\in H_{\alpha}^{+}(\Omega)$ is nonnegative according to  \cite[Theorem 3.1]{Sh}. Let
$E$ be a set as in Theorem \ref{Theorem 1.5} (ii)
and $$
E_{\varepsilon}=\{x:f(x)\geq \overline{f(x_{0})}-\varepsilon\}\cap \{\Omega\backslash E\}
$$ for each $\varepsilon>0.$
Then $C_{\alpha,p}(B\cap E_{\varepsilon},2B)>0$ for each ball $B=B(x_{0},r)$ since $E_{\varepsilon}$ is not $(\alpha,p)-$thin at $x_{0}$. Accordingly, $u(x_{0})\geq \overline{f(x_{0})}-\varepsilon$, and hence $u(x_{0})\geq \overline{f(x_{0})}$.

On the other hand, \cite[Theorem 1.1 \& Remark 4.2]{CKP2} gives
$$
\esss _{B(x_{0},r)}(u-\xi)^{+}\leq \frac{C(u-\xi)^{+}}{|B(x_{0},2r)|}dx=\frac{C(u-\min\{u,\xi\})}{|B(x_{0},2r)|}dx\longrightarrow  0\quad\hbox{as}\quad r\longrightarrow  0
$$
for fixed $\xi>u(x_{0})$ and $r>0$ small enough, where the last approximation follows from \cite[Theorem 13]{KKP4}. On account of the above analysis and \cite[Corollary 3.9]{Sh}, we have $\limsup_{x\longrightarrow  x_{0}}u(x)\leq \xi$, which gives the desired result by letting $\xi\longrightarrow  u(x_{0}).$

Next, we proceed to show (ii). Assume that there exists $\varepsilon>0$ such that $F_{\varepsilon}$ is $(\alpha,p)-$thin at $x_{0}$, then Theorem \ref{Theorem 1.1} gives the existence of a bounded function $u\in H_{\alpha}^{+}(B(x_{0},r_{0}))$ with $u=\sup f$ on $F_{\varepsilon}\backslash \{x_{0}\}$ and $u(x_{0})=\overline{f(x_{0})}-\frac{\varepsilon}{2}$. By \cite[Theorem 15]{P}, we may assume that $u\in W^{\alpha,p}(B)$ and $u>\overline{f(x_{0})}-\varepsilon$ in $B$.

Let $v\in H_{\alpha}^{+}(\Omega)$ be the solution to $\Phi^{x_0}_{f,u_{0}}$ with $\Omega=B$ such that $u-v\in W^{\alpha,p}_{0}(B)$. Then \cite{KKP4} implies that $v\leq u$
in $B$ since $u\geq f$ q.e.. Accordingly, we get
$$
\liminf_{x\longrightarrow  x_{0}}v(x)\leq \essliminf_{x\longrightarrow  x_{0}}u(x)=u(x_{0})<\overline{f}(x_{0})=\inf_{r>0}(\alpha,p)-\esss_{B(x_{0},r)} f\leq \limsup_{x\longrightarrow  x_{0}}v(x)
$$
for $v\geq f$ q.e., which shows that $v$ can not be continuous at $x_{0}$.

\section{Beyond}\label{s8}

We end this paper with more regulairty conditions: the cotinuitity of fractional superharmonic functions,
the fractional resolutivity, a connection between  $(\alpha,p)-$potentials and  $(\alpha,p)-$Perron solutions by using Lemmas \ref{lemma 3.4} \& \ref{lemma 3.6},
and the existence of a capacitary function for an arbitrary condenser by applying Theorem \ref{theorem 3.2} and Lemma \ref{lemma 5.1}, which can be viewed as a general version of \cite[Theorem 1.2 (ii)]{SX3}.

\subsection{The continuity of fractional superharmonic functions}\label{s6}

With the help of a generalized comparison lemma, we can deduce the continuity of $(\alpha,p)-$superharmonic functions from Theorem \ref{Theorem 1.1}.

\begin{theorem}\label{Theorem 1.5}
Let $u\in H_{\alpha}^{+}(\Omega)$. Then the following three statements hold.
\begin{itemize}	
\item[\rm (i)] $u$ is real-valued and continuous at $x_{0}$ if and only if for each $\varepsilon>0$, there is an $r>0$ such that
$$
\textbf{W}_{\alpha,p}^{\mu}(x,r)<\frac{\varepsilon}{ C}-\hbox{Tail}(u^{-};x,r),
$$
where $$\textbf{W}_{\alpha,p}^{\mu}(x,r)=:\int_{0}^{r}\left(\frac{\mu(B(x,t))}{t^{n-\alpha p}}\right)^{\frac{1}{p-1}}\frac{dt}{t}$$ is the fractional Wolff potential for a Radon measure $\mu$ (\cite{KMS}),
and
$$\hbox{Tail}(u^{-};x,r)=:\left(r^{\alpha p}\int_{\mathbb{R}^n\backslash B(x,r)}|u^{-}(x_0)|^{p-1}|x-x_0|^{-n-\alpha p}dx_0\right)^{\frac{1}{p-1}}$$
in the nonlocal tail of function $u^{-}$ in the ball of radius $r>0$ centred in $x\in \mathbb{R}^n.$

\item[\rm (ii)] There exists a set $E$ such that $E$ is  $(\alpha,p)-$thin at $x_{0}\in \Omega$ and that $u|_{\Omega\backslash E}$ is continuous at $x_{0}$;

\item[\rm (iii)] If a finite point $x_{0}\in \partial\Omega$ is $(\alpha,p)-$ regular with $C_{\alpha, p}(\{x_{o}\}, \Omega)=0$ and $u\in H_{\alpha}^{+}(U(x_{0}))$, then $\lim_{x\longrightarrow  x_{0}\ \& \ x\in \Omega^{c}}u(x)=u(x_{0}).$
\end{itemize}
\end{theorem}

We begin with the generalized comparison lemma which is critical to the proof of Theorem \ref{Theorem 1.5} (iii).
\begin{lemma}\label{lemma 5.1}
Let
$$
\begin{cases}
\Omega\subset \mathbb{R}^{n}\,\,\hbox{be bounded};\\
E\subset \partial \Omega\,\,\hbox{and}\,\, C_{\alpha,p}(E,\Omega)=0;\\
(u,v)\in (H_{\alpha}^{+}(\Omega)\times H_{\alpha}^{-}(\Omega))\cap (\dot{W}^{\alpha,p}(\Omega)\times \dot{W}^{\alpha,p}(\Omega)).
\end{cases}
$$
Then $v\leq u$ in $\Omega$.
\end{lemma}

\begin{proof}
Without loss of generality, we  assume that $E$ is compact and $u\in \dot{W}^{\alpha,p}(\Omega)$ since the set $$\left\{x\in \partial\Omega: \liminf_{y\longrightarrow  x}u(y)+\varepsilon>\limsup_{y\longrightarrow  x}v(y)\right\}$$ is open on $\partial\Omega$ for $\varepsilon>0$. We claim that there exists a decreasing sequence $\{\varphi_{i}\}\subset C(\Omega\cap E)\cap W^{\alpha,p}(\Omega)$ such that
$$
\begin{cases}
\varphi_{i}\in [0,M]\quad\hbox{for}\quad M=\sup |u|+\sup |v|;\\
\varphi_{i}=M\quad\hbox{on}\quad E;\\
\|\varphi_{i}\|_{W^{\alpha,p}(\Omega)}\longrightarrow  0\quad\hbox{as}\quad i\longrightarrow  \infty.
\end{cases}
$$
This can be done by picking a nonnegative sequence $\{\eta_{i}\}\subset C(\Omega\cap E)\cap W^{\alpha,p}(\Omega)$ with $\eta_{i}=M$ on $E$ and $\eta_{i}\longrightarrow  0$
in $W^{\alpha,p}(\Omega)$. In fact, we can pick $\{\varphi_{i}\}$ as follows
$$
\begin{cases}
\varphi_{1}=\min\{M, \eta_{1}\};\\
\varphi_{i+1}=\min\{\varphi_{i}, \eta_{j}\}\quad\hbox{with}\quad j\quad\hbox{is choosen large enough to gurantee}\quad\|\varphi_{i+1}\|_{W^{\alpha,p}(\Omega)}\leq\frac{1}{2}\|\varphi_{i}\|_{W^{\alpha,p}(\Omega)}
\end{cases}
$$
since $\min\{\varphi_{i}, \eta_{j}\}\longrightarrow  0$ in $W^{\alpha,p}(\Omega)$ as $j\longrightarrow  \infty$ according to  \cite[Lemma 2.1]{SX3}. Let
$$
\begin{cases}
\xi_{i}=u+\varphi_{i};\\
u_{i}\quad\hbox{be the solution to}\quad\Phi_{\psi_{i},\psi_{i}(\Omega)}.
\end{cases}
$$
It follows from \cite[Theorem 13]{KKP4} that $u_{i}\geq \xi_{i}$ in $\Omega$. Hence, we get $\limsup_{y\longrightarrow  x}v(y)\leq \liminf_{y\longrightarrow  x}u_{i}(y)$ for $x\in \partial\Omega$, which shows $u_{i}\geq v$ according to  \cite[Proposition 2.20]{Sh}.
On the other hand,  \cite[Proposition 2.15]{Sh} implies $u_{i}\longrightarrow  u.$  A further application of \cite[Theorem 13]{KKP4} gives $u\geq v$ in $\Omega$ as desired.
\end{proof}

We now turn to the proof of Theorem \ref{Theorem 1.5}.
The sufficiency of (i) is \cite[Theorem 1.5]{KMS}.
We only need to prove the necessity.  Without  loss of generality, we assume $u(x_{0})=0$ for $u(x_{0})<\infty$ according to  \cite[Theorem 1.2 \& Theorem 1.3]{KMS} and choose
$r_{0}>0$ with $u(x)>-\varepsilon$ for $x\in B(x_{0}, 4r_{0})$ since $u$ is l.s.c. Then, we get,  for $x\in B(x_{0},r)$ with $r<r_{0}$,
$$
u(x)<C\inf_{B(x,r)}u+C\textbf{W}_{\alpha,p}^{\mu}(x,r)+C\hbox{Tail}(u^{-};x_{0},r)\leq C\varepsilon
$$
as desired.

(ii) can be proved by applying Theorem \ref{Theorem 1.1}, which will be divided into two cases.

Case 1. $u(x_{0})=\infty$. It is easy to obtain the desired result by choosing $E=\phi$.

Case 2. $u(x_{0})<\infty$. Denote by $E_{i}=\left\{x\in \Omega: u(x)-u(x_{0})>\frac{1}{i}\right\}$ for $i\in \mathbb{N}^{+}$. Then $E_{i}$ is $(\alpha,p)-$thin at $x_{0}$ by Theorem \ref{Theorem 1.1}. Using \cite[Proposition 3.1]{Me}, we can show  that there exists a decreasing sequence $r_{i}\longrightarrow  0$ such that the set $E=\cup_{i}(E_{i}\cap B(x_{0},r_{i}))$ is $(\alpha,p)-$thin at $x_{0}$, which  implies our desired result since $u$ is l.s.c.

Now, we prove (iii) by a contradiction. If (iii) is not true, assume that there exists a function $u\in H_{\alpha}^{+}(U(x_{0}))$ such that $\lim_{x\longrightarrow  x_{0}, x\in \Omega^{c}}u(x)>u(x_{0})$.
For some ball $B$ centered at $x_{0}$, we may assume that $u\in W^{\alpha, p}(3B)$ and  $u=1$ in $\Omega^{c}\cap (\overline{B} \backslash\{x_{0}\})$. Lemma \ref{lemma 5.1} now shows that $\overline{H}_{f}^{\alpha}\leq u$ in $2B\backslash (\overline{B}\cap \Omega^{c})$ by picking
$$
\begin{cases}
f\in C^{\infty}(\mathbb{R}^{n})\quad\hbox{with}\quad f=1\quad\hbox{in}\quad\overline{B};\\
f\leq u\quad\hbox{on}\quad\partial (2B).
\end{cases}
$$
It follows from \cite[Theorem 13]{KKP4} that
$$
1=\lim_{x\longrightarrow  x_{0}, x\in \Omega}\overline{H}_{f}^{\alpha}(x)\leq \liminf_{x\longrightarrow  x_{0}, x\in \Omega}u(x)=u(x_{0})<1,
$$
which is a contradiction.

\subsection{The fractional resolutivity}\label{s7}

A function $f:\partial\Omega\longrightarrow  [-\infty,+\infty]$ is called $(\alpha,p)$-resolutive if $\overline{H}_{f}^{\alpha}= \underline{H}_{f}^{\alpha}\in H_{\alpha}(\Omega)$. \cite[Theorem 2]{KKP4} shows that $f$ is $(\alpha,p)$-resolutive if $f$ is continuous. The resolutivity of  $f$ does not imply the Perron regularity of a boundary point. But, the inverse is true. More generally, we have the following result.
\begin{theorem}\label{Theorem 1.6}
Let $f:\partial\Omega\longrightarrow  [-\infty,+\infty]$. Then $f$ is $(\alpha,p)-$resolutive if one of the following conditions holds.
\begin{itemize}	
\item[\rm (i)] There is a bounded $u\in H_{\alpha}(\Omega)$ with $\lim_{x\longrightarrow  x_{0}}u(x)=f(x_{0})$ for each $x_{0}\in \partial\Omega$;

\item[\rm (ii)] $C_{\alpha,p}(\Omega^{c})>0$ and $f\in C(\partial\Omega)$;

\item[\rm (iii)] $\Omega$ is Perron regular and $f\in C(\partial\Omega)$;

\item[\rm (iv)] $\Omega$ is Perron regular, $f$ is bounded and lower semicontinuous on $\partial\Omega$.
\end{itemize}
\end{theorem}

In order to
prove Theorem \ref{Theorem 1.6}, we  need to establish the following uniform convergence lemma for upper $(\alpha,p)$-Perron solutions.

\subsubsection{The uniform convergence of upper $(\alpha,p)$-Perron solutions}
We call a set of functions  $\mathcal{E}$ is downward directed if for functions $g_{1}, g_{2}\in \mathcal{E}$, there exists a function $g\in \mathcal{E}$ such that $g\leq \min\{g_{1}, g_{2}\}.$

\begin{lemma}\label{lemma 6.1}
\begin{itemize}
\item[\rm (i)]  Let $G$ be a downward directed family of upper semicontinuous functions $g: \partial\Omega\longrightarrow  [-\infty, +\infty]$ and let $f=\inf G$. Then
$
\overline{H}_{f}^{\alpha}=\inf\left\{\overline{H}_{g}^{\alpha}: g\in G\right\}.
$

\item[\rm (ii)] Suppose that $f_{i}: \partial\Omega\longrightarrow  [-\infty, +\infty]$ is a decreasing sequence of upper semicontinuous functions and $\lim_{i}f_{i}=f$. Then
$\overline{H}_{f}^{\alpha}=\lim_{i\longrightarrow  \infty}\overline{H}_{f_{i}}^{\alpha}$.

\item[\rm (iii)]
Let
$
\{f_{i}\}: \partial \Omega\longrightarrow  \mathbb{R}$ be resolutive and
$f_{i}\longrightarrow  f$ uniformly.
Then $f$ is $(\alpha,p)-$resolutive and $\overline{H}_{f_{i}}^{\alpha}\longrightarrow  \overline{H}_{f}^{\alpha}$ as $i\longrightarrow  \infty$.
\end{itemize}
\end{lemma}

\begin{proof}
(i). For abbreviation, we denote $h=\inf_{G}\overline{H}_{g}^{\alpha}$. Next, we  show $\overline{H}_{f}^{\alpha}=h.$ The inequality $\overline{H}_{f}^{\alpha}\leq h$ in $\Omega$ is straightforward. To show the reverse inequality, let $u\in \mathcal{U}_{f}^{\alpha}.$ Then,  for fixed $\varepsilon>0$, the set $$\mathcal{F}=\left\{y\in \partial\Omega: \liminf_{x\longrightarrow  y}u(x)+\varepsilon>g(y), g\in G\right\}$$
is open and cover $\partial \Omega$ according to  the upper semicontinuity of $G$. The fact that $\partial \Omega$ is compact in $\overline{\mathbb{R}^{n}}$ and $G$ is downward directed allow us to find a function $g\in G$ with
$\lim_{x\longrightarrow  y}u(x)+\varepsilon>g(y)\quad\forall y\in \partial\Omega.$
Namely, $u+\varepsilon\in \mathcal{U}_{g}^{\alpha}$, and hence $u+\varepsilon\geq  \overline{H}_{g}^{\alpha}\geq h.$ Finally we derive $\overline{H}_{f}^{\alpha}+\varepsilon\geq h$, which  shows that $\overline{H}_{f}^{\alpha}\geq h$ by the arbitrariness of $\varepsilon$ and the desired result is obtained.

(ii) follows from (i).

(iii). Since $|f_{i}-f|<\varepsilon$ for given $\varepsilon>0$ and $i$ large enough,
$
\overline{H}_{f}^{\alpha}-\varepsilon\leq \overline{H}_{f_{i}}^{\alpha}=\underline{H}_{f_{i}}^{\alpha}\leq \underline{H}_{f}^{\alpha}+\varepsilon.
$
Consequently, $\lim_{i\longrightarrow  \infty}\overline{H}_{f_{i}}^{\alpha}=\underline{H}_{f}^{\alpha}=\overline{H}_{f}^{\alpha}$.
Furthermore, we get  $\overline{H}_{f}^{\alpha}\in H_{\alpha}(\partial\Omega)$ since $\overline{H}_{f}^{\alpha}$ is finite.
\end{proof}

\subsubsection{ Proof of Theorem \ref{Theorem 1.6}}.
(i) can be derived immediately from the comparison principle, see for example \cite[Proposition 2.20]{Sh}.

(ii) follows from  (iii) of Lemma \ref{lemma 6.1}. Indeed, without loss of generality, we assume  $f\in C(\overline{\mathbb{R}^{n}})$ and $f(\infty)=0$  according to the Tietze extension theorem. Then,  we get the desired result since $f$ can be approximated uniformly by $\varphi_{j}\in C_{c}^{\infty}(\mathbb{R}^{n})$ and $\varphi_{j}$ are resolutive according to   (i) of Lemma \ref{lemma 3.6}.

(iii) is  a by-product of Lemma \ref{lemma 3.2} (iii) and Theorem \ref{Theorem 1.4}.

(iv) will be proved by showing $\underline{H}_{f}^{\alpha}\geq \overline{H}_{f}^{\alpha}$. Let $\{f_{i}\}\subset C(\partial \Omega)$ be increasing such that $f_{i}\longrightarrow  f$ on $\partial \Omega.$ Then,  we have
$$
\liminf_{x\longrightarrow  y}\underline{H}_{f}^{\alpha}(x)\geq \lim_{x\longrightarrow  y}\underline{H}_{f_{i}}^{\alpha}(x)=f_{i}(y)\quad\forall y\in \partial \Omega.
$$
Accordingly, we get $\underline{H}_{f}^{\alpha}\geq \overline{H}_{f}^{\alpha}$ by the fact $\underline{H}_{f}^{\alpha}\in \mathcal{U}_{f}^{\alpha}$ since $f_{i}\longrightarrow  f.$

\subsection{A connection between  $(\alpha,p)-$potentials and  $(\alpha,p)-$Perron solutions}

\begin{theorem}\label{theorem 3.2}
Let $u\in H_{\alpha}^{+}(\Omega)$ be nonnegative and $\Omega\subset \mathbb{R}^{n}$ be open. Then the following statements hold.
\begin{itemize}
\item[\rm (i)]  $\widehat{\mathcal{B}}^{u}_{E}(\Omega)=\overline{H}_{f}^{\alpha}(\Omega\backslash E)$ in $\Omega\backslash E$ with $E\subset \Omega$ is relatively closed and
$
f=\begin{cases}u\quad\hbox{on}\quad\partial E\cap \Omega;\\
0\quad\hbox{on}\quad\partial \Omega.
\end{cases}
$
\item[\rm (ii)]  $$
\begin{cases}\overline{H}_{f}^{\alpha}-f\in \dot{W}^{\alpha,p}_{0}(\Omega\backslash E)\quad\hbox{if}\, \,  f\in C(\overline{\Omega})\cap \dot{W}^{\alpha,p}(\Omega);\\
\overline{H}_{f}^{\alpha}-f\in W^{\alpha,p}_{0}(\Omega\backslash E)\quad\hbox{if}\, \, \Omega\,\, \hbox{is bounded and}\,  f\in C(\overline{\Omega})\cap \dot{W}^{\alpha,p}(\Omega).
\end{cases}
$$

\item[\rm (iii)] $\lim_{x\longrightarrow  x_{0}}\widehat{\mathcal{B}}^{u}_{E}(x)=0$ for $E\Subset \Omega$ and $x_{0}\in \partial\Omega$ is Perron regular. Particularly, $\lim_{x\longrightarrow  x_{0}}\widehat{\mathcal{B}}^{u}_{E}(x)=0$ q.e. on $\partial \Omega$.
\end{itemize}
\end{theorem}
\begin{proof}
(i). For the first assertion, it remains to prove that $\overline{H}_{f}^{\alpha}\geq \widehat{\mathcal{B}}^{u}_{E}$ in $\Omega\backslash E$ since the reverse inequality is obviously. To do so, by choosing $v\in \mathcal{U}_{f}^{\alpha}$ and writing
$$
w=\begin{cases}\min\{u,v\}\quad\hbox{in}\quad\Omega\backslash E;\\
u\quad\hbox{on}\quad E,
\end{cases}
$$
we get $w\in H_{\alpha}^{+}(\Omega)$ according to  \cite[Proposition 2.2]{Sh}.  Hence, we obtain $w\geq \widehat{\mathcal{B}}^{u}_{E},$ and finally $v\geq \widehat{\mathcal{B}}^{u}_{E}$ as desired.

The second assertion follows from (i) of Lemma \ref{lemma 3.4}  and   (ii) of Lemma \ref{lemma 3.6}.

(iii). By choosing a polyhedron $O\Subset \Omega$ with $E\Subset O$ and applying the Poisson modification in a neighborhood of $\partial O$, there exists a function $v\in \Psi_{E}^{u}(\Omega)$ with $v$ bounded on $\partial O$. Therefore, we  have
$0\leq \widehat{\mathcal{B}}^{u}_{E}\leq \widehat{\mathcal{B}}^{v}_{\overline{O}},$  which implies the desired result since $\lim_{x\longrightarrow  y}\widehat{\mathcal{B}}^{v}_{\overline{O}}(x)=0$ for any $y\in \partial\Omega.$ Moreover, for $y\in \partial\Omega,$ (i) implies that $y$ is Perron regular. Furthermore,   it follows from  \cite[Proposition 3.2]{SX3} that   $\lim_{x\longrightarrow  x_{0}}\widehat{\mathcal{B}}^{u}_{E}(x)=0$ q.e. on $\partial \Omega.$
\end{proof}

\subsection{The existence of a capacitary function for an arbitrary condenser}
\begin{theorem}\label{theorem 3.3}
Let $E\subset \Omega$  with $ C_{\alpha,p}(E, \Omega)<\infty$  and
$u=\widehat{\mathcal{B}}^{1}_{E}(\Omega).$
Then
$$
\begin{cases}u\in \dot{W}^{\alpha,p}(\Omega)\quad\hbox{and}\quad C_{\alpha,p}(E, \Omega)=[u]_{\dot{W}^{\alpha,p}(\Omega)};\\
u\in W^{\alpha,p}(\Omega)\quad\hbox{if}\quad\Omega\quad\hbox{is bounded}.
\end{cases}
$$
\end{theorem}
\begin{proof}
The proof can be divided into three steps.

Step 1-showing $[v_{i}]_{\dot{W}^{\alpha,p}(\Omega)}=C_{\alpha,p}(\overline{O}_i, \Omega)$, where
$$
\begin{cases}O\subset \Omega\quad\hbox{is open with}\quad C_{\alpha,p}(O, \Omega)<\infty;\\
\{O_{i}\}\quad\hbox{is a sequence of  polyhedra with}\quad O_{1}\Subset O_{2}\Subset \cdots \Subset O;\\
v_{i}=\widehat{\mathcal{B}}^{1}_{\overline{O_{i}}}(\Omega)\quad i=1,2,3,\cdots .
\end{cases}
$$
In fact, we deduce from Theorem \ref{theorem 3.2} (i) that $v_{i}-\varphi\in \dot{W}^{\alpha,p}_{0}(\Omega\backslash \overline{O_{i}})$ if $\varphi$ is admissible for the condenser $(\overline{O_{i}}, \Omega)$, namely, $\varphi\in X_{0}^{\alpha,p}(\overline{O_{i}}, \Omega)$. Since $\varphi_{i}\in H_{\alpha}(\Omega\backslash \overline{O_{i}})$, we get
$$
C_{\alpha,p}(\overline{O_{i}}, \Omega)\leq [v_{i}]_{\dot{W}^{\alpha,p}(\Omega)}\leq [\varphi]_{\dot{W}^{\alpha,p}(\Omega)}.
$$
Then,  the desired result follows by taking the infimum over all $\varphi$.

Step 2-proving $[u]_{\dot{W}^{\alpha,p}(\Omega)}=C_{\alpha,p}(E, \Omega)$ if $E=O$ is open. Using Lemma \ref{lemma 4.1}, we derive   that $v_{i}\longrightarrow  \widehat{\mathcal{B}}^{1}_{E}$ weakly in $\dot{W}^{\alpha,p}(\Omega)$ since
$$
[v_{i}]_{\dot{W}^{\alpha,p}(\Omega)}=C_{\alpha,p}(\overline{O_{i}}, \Omega)\leq C_{\alpha,p}(O, \Omega)<\infty.
$$
Then $\widehat{\mathcal{B}}^{1}_{E}$ can be approximated in $\dot{W}^{\alpha,p}(\Omega)$ by functions $\varphi\in X_{0}^{\alpha,p}(\overline{O_{i}}, \Omega)$ for each $i$ by the Mazar lemma. Consequently, \cite[(5.25)]{HKM2} gives
$$
C_{\alpha,p}(O, \Omega)=\lim_{i\longrightarrow \infty}C_{\alpha,p}(\overline{O_{i}}, \Omega)\leq[u]_{\dot{W}^{\alpha,p}(\Omega)}\leq
\liminf_{i\longrightarrow  \infty}[v_{i}]_{\dot{W}^{\alpha,p}(\Omega)}=\lim_{i\longrightarrow \infty}C_{\alpha,p}(\overline{O_{i}}, \Omega)= C_{\alpha,p}(O, \Omega)
$$
as desired.

Step 3-checking $[u]_{\dot{W}^{\alpha,p}(\Omega)}=C_{\alpha,p}(E, \Omega)$ for a general set $E\subset\Omega$. Let
$$
\begin{cases}u_{i}\in C_{c}^{\infty}(\Omega)\quad\hbox{be increasing and nonnegative};\\
u_{i}\longrightarrow  u\quad\hbox{in}\quad\Omega;\\
E\subset O\subset \Omega\quad\hbox{for a nonempty open set}\quad O\quad\hbox{with}\quad C_{\alpha,p}(O, \Omega)<\infty.
\end{cases}
$$
Then we may assume that $u_{i}<\widehat{\mathcal{B}}^{1}_{O}$ in $\Omega$ by replacing $u_{i}$ with $(1-\delta_{i})u_{i}$ if necessary for $\delta_{i}\longrightarrow  0$
decreasingly. Exhausting $\Omega$ by polyhedra $U_{1}\Subset U_{2}\Subset \cdots\Subset \Omega$, there is an index $i_{j}$ with $\widehat{\mathcal{B}}^{1}_{O\cap U_{i_{j}}}(\Omega)>u_{i}$ since $\widehat{\mathcal{B}}^{1}_{O\cap U_{i}}(\Omega)$ increases to $\widehat{\mathcal{B}}^{1}_{O}(\Omega)$ after the observation
$\widehat{\mathcal{B}}^{1}_{O\cap U_{i}}(\Omega)>0$ for $i$ large enough.
Likewise, setting $O_{j}=O\cap U_{i_{j}}$, there is a $k_{j}$ with $\widehat{\mathcal{B}}^{1}_{O_{j}}(U_{k})\geq u_{j}$ in $U_{k}$ for $k\geq k_{j}$ since $\widehat{\mathcal{B}}^{1}_{O_{j}}(U_{k})$ increases to $\widehat{\mathcal{B}}^{1}_{O_{j}}(\Omega)$.

By setting
$$
\begin{cases}v_{j,k}=\widehat{\mathcal{B}}^{1}_{O_{j}}(U_{k});\\
u_{j}=\widehat{\mathcal{B}}^{\varphi_{j}}(\Omega);\\
u_{j,k}=\widehat{\mathcal{B}}^{\varphi_{j}}(U_{k}),
\end{cases}
$$
we obtain that $v_{j,k}$ increases to $\widehat{\mathcal{B}}^{1}_{O_{j}}(\Omega)$ as $k\longrightarrow \infty$. Then by extending $v_{j,k}$ as zero to $\Omega\backslash U_{k} $, one has
$$
[v_{j,k}]_{\dot{W}^{\alpha,p}(\Omega)}=C_{\alpha,p}(O_{j}, U_{k})\longrightarrow  C_{\alpha,p}(O_{j}, \Omega)<\infty,
$$
and hence  $v_{j,k}\longrightarrow  \widehat{\mathcal{B}}^{1}_{O_{j}}(\Omega)$ weakly in $\dot{W}^{\alpha,p}(\Omega)$ by Lemma \ref{lemma 4.1}.  Finally, we get that
$\lim_{k\longrightarrow  \infty}u_{j,k}=u_{j}$. Moreover, the quasi minimizing property shows that
$
[u_{j,k}]_{\dot{W}^{\alpha,p}(\Omega)}\leq [v_{j,k}]_{\dot{W}^{\alpha,p}(\Omega)},
$
and hence
$$
(v_{j,k}-u_{j,k})\longrightarrow  (\widehat{\mathcal{B}}^{1}_{O_{j}}(\Omega)-u_{j})\quad\hbox{weakly in}\quad\dot{W}^{\alpha,p}(\Omega).
$$
Thus, we get
$$
0\leq [u_{j}(v_{j,k}-u_{j,k})]_{\dot{W}^{\alpha,p}(\Omega)}\longrightarrow  [u_{j}(\widehat{\mathcal{B}}^{1}_{O_{j}}(\Omega)-u_{j})]_{\dot{W}^{\alpha,p}(\Omega)}\quad\hbox{as}\quad k\longrightarrow \infty.
$$
Accordingly, we reach
$$
[u_{j}]_{\dot{W}^{\alpha,p}(\Omega)}\leq [\widehat{\mathcal{B}}^{1}_{O_{j}}]_{\dot{W}^{\alpha,p}(\Omega)}=C_{\alpha,p}(O_{j}, \Omega)\leq C_{\alpha,p}(O, \Omega)<\infty.
$$
A further application  of Lemma \ref{lemma 4.1} shows that $u_{j}\longrightarrow  u$ weakly in $\dot{W}^{\alpha,p}(\Omega)$ and $u\in \dot{W}^{\alpha,p}_{0}(\Omega)$. Furthermore,  $u\in W^{\alpha,p}_{0}(\Omega)$ if $\Omega$ is bounded according to  \cite[Lemma 2.2]{Sh}. The weak lower semicontinuity of norms gives that
$$
[u]_{\dot{W}^{\alpha,p}(\Omega)}\leq \liminf_{j\longrightarrow \infty}[u_{j}]_{\dot{W}^{\alpha,p}(\Omega)}\leq \lim_{j\longrightarrow \infty}C_{\alpha,p}(O_{j}, \Omega)= C_{\alpha,p}(O, \Omega).
$$
Consequently, we obtain  $[u]_{\dot{W}^{\alpha,p}(\Omega)}\leq C_{\alpha,p}(E, \Omega)$ by taking the infimum over all open neighborhoods $O$ of $E$ in $\Omega$.

What is left is to show that $[u]_{\dot{W}^{\alpha,p}(\Omega)}\geq C_{\alpha,p}(E, \Omega)$. For $0<\varepsilon<1,$ denote by
$$ O_{\varepsilon}=\{u>1-\varepsilon\} \quad \hbox{and}\quad  E^{'}=\{x\in E, u(x)=1\}.$$
Then $O_{\varepsilon}\supset E^{'}$ is open and $ C_{\alpha,p}(E', \Omega)= C_{\alpha,p}(E, \Omega)$ by \cite[Theorem 1.3]{SX3}. It follows from  step 1 that $[\widehat{\mathcal{B}}^{1}_{\overline{O}}]_{\dot{W}^{\alpha,p}(\Omega)}=C_{\alpha,p}(\overline{O}, \Omega)$ for $O\Subset O_{\varepsilon}$,  which is a polyhedron. On the other hand, it can be deduced from the fact $u\in \dot{W}^{\alpha,p}_{0}(\Omega)$ and $\widehat{\mathcal{B}}^{1}_{\overline{O}}\in \dot{W}^{\alpha,p}_{0}(\Omega)$ that $\left(\min\left\{1, \frac{u}{1-\varepsilon}\right\}-\widehat{\mathcal{B}}^{1}_{\overline{O}}\right)\in \dot{W}^{\alpha,p}_{0}(\Omega\backslash \overline{O})$. Accordingly, we derive
$$
C_{\alpha,p}(\overline{O}, \Omega)=[\widehat{\mathcal{B}}^{1}_{\overline{O}}]_{\dot{W}^{\alpha,p}(\Omega)}\leq (1-\varepsilon)^{-p}[u]_{\dot{W}^{\alpha,p}(\Omega)}
$$
since $\widehat{\mathcal{B}}^{1}_{\overline{O}}\in H_{\alpha}(\Omega\backslash \overline{O})$.

Finally,  we obtain $C_{\alpha,p}(E, \Omega)\leq C_{\alpha,p}(O_{\varepsilon}, \Omega)\leq (1-\varepsilon)^{-p}[u]_{\dot{W}^{\alpha,p}(\Omega)}$ by taking the supremum over all polyhedra $O$ in $O_{\varepsilon}$, which implies the desired inequality by taking $\varepsilon\longrightarrow  0$.
\end{proof}



\end{document}